\documentclass[10pt]{article}
\usepackage{amsthm,amsfonts,amsmath,amscd,amssymb,epsfig,psfrag}
\usepackage{enumerate}

\theoremstyle{plain}

\newtheorem{cor}{Corollary}
\newtheorem{lemma}{Lemma}
\newtheorem{fact}{Fact}
\newtheorem{problem}{Problem}
\newtheorem{question}{Question}
\newtheorem{prop}{Proposition}
\newtheorem{example}{Example}
\newtheorem{conjecture}{Conjecture}

\theoremstyle{remark}
\newtheorem*{remark}{Remark}

\theoremstyle{definition}

\newcommand{\ol}{\overline}

\newcommand{\p}{\partial}

\newcommand{\om}{\omega}
\newcommand{\Om}{\Omega}
\newcommand{\eps}{\varepsilon}
\newcommand{\into}{\hookrightarrow}

\newcommand{\N}{{\mathbb{N}}}
\newcommand{\Z}{{\mathbb{Z}}}
\newcommand{\R}{{\mathbb{R}}}
\newcommand{\C}{{\mathbb{C}}}
\newcommand{\Q}{{\mathbb{Q}}}
\renewcommand{\P}{{\mathbb{P}}}

\newcommand{\ind}{{\rm ind}}

\newcommand{\vol}{{\rm vol}}

\newcommand{\const}{{\rm const}}

\renewcommand{\min}{{\rm min}}

\newcommand{\EH}{{\rm EH}}
\newcommand{\HZ}{{\rm HZ}}

\newcommand{\graph}{{\rm graph}}

\newcommand{\CC}{\mathcal{C}}
\newcommand{\LL}{\mathcal{L}}
\newcommand{\FF}{\mathcal{F}}

\newcommand{\HH}{\mathcal{H}}

\renewcommand{\AA}{\mathcal{A}}
\newcommand{\GG}{\mathcal{G}}

\newcommand{\cs}{\mathcal{S}}

\newcommand{\Symp}{\it{Symp}}
\newcommand{\Op}{\it{Op}}
\newcommand{\Conv}{\it{Conv}}
\newcommand{\Ell}{\it{Ell}}
\newcommand{\LinEll}{\it{LinEll}}
\newcommand{\Pol}{\it{Pol}}
\newcommand{\LinPol}{\it{LinPol}}
\newcommand{\RCT}{\it{RCT}}
\newcommand{\proofend}{\hspace*{\fill} $\Box$\\}
\newcommand{\diam}{\hspace*{\fill} $\Diamond$}
\newcommand{\comment}[1]{}
\title{Quantitative symplectic geometry}
\author{K.~Cieliebak, H.~Hofer, J.~Latschev and F.~Schlenk
\footnote{The research of the first author was partially supported
by the DFG grant Ci 45/2-1. The research of the second  author was
partially supported by the NSF Grant DMS-0102298. The third author
held a position financed by the DFG grant Mo 843/2-1. The fourth
author held a position financed by the DFG grant Schw 892/2-1. }}
\date{\today}
\begin{document}
\maketitle


A symplectic manifold $(M,\omega)$ is a smooth manifold $M$
endowed with a non-degenerate and closed $2$-form $\omega$.
By Darboux's Theorem such a manifold looks locally like an open set in
some $\R^{2n} \cong \C^{n}$ with the standard symplectic form
\begin{equation}  \label{def:omega0}
\omega_0 \,=\, \sum_{j=1}^n dx_j \wedge dy_j ,
\end{equation}
and so symplectic manifolds have no local invariants.
This is in sharp contrast to Riemannian manifolds, for which the
Riemannian metric admits various curvature invariants. Symplectic
manifolds do however admit many global numerical invariants, and
prominent among them are the so-called symplectic capacities.

Symplectic capacities were introduced in 1990 by I.~Ekeland and
H.~Hofer~\cite{EH-90a,EH-90b} (although the first capacity was in fact
constructed by M.~Gromov~\cite{Gr}). Since then, lots of new capacities
have been defined
\cite{CM-03a,FHW94,FS,Hofer90,HZ90,LP2004,Lu98,Schwarz00,V}
and they were further studied in
\cite{Bates95,Bates98,BC,BPS,CM-03b,EH93,EM92,FHV90,FGS,G03,GG,GK99,He00,
He-inner,Hermann,Hofer93,HV92,HZ,
Jiang93,LM95,LM95b,LP,Lu02,Mac,Mac2,Mac3,MaS,MMT,Mc91b,MS01,McTr,Neduv01,Sch,Sch1,Siburg93,Sikorav90,Tokieda,Viterbo89,Viterbo90}.
Surveys on symplectic capacities are
\cite{Hofer90b,HZ,Lalonde97,Mc99,Viterbo89}.
Different capacities are defined in different ways, and so relations
between capacities often lead to surprising relations between different
aspects of symplectic geometry and Hamiltonian dynamics. This is
illustrated in \S~\ref{s:caps}, where we discuss some examples
of symplectic capacities and describe a few consequences of their existence.
In \S~\ref{s:relations} we present an attempt to better understand
the space of all symplectic capacities, and discuss some further general
properties of symplectic capacities.
In \S~\ref{s:EP}, we describe several new relations between
certain symplectic capacities on ellipsoids and polydiscs.
Throughout the discussion we mention many open problems.

As illustrated below, many of the quantitative
aspects of symplectic geometry can be formulated in terms of
symplectic capacities. Of course there are other numerical invariants
of symplectic manifolds which could be included in a discussion of
quantitative symplectic geometry, such as the invariants derived from
Hofer's bi-invariant metric on the group of Hamiltonian diffeomorphisms,
\cite{Hofer90,Polterovich-spec,Polterovich-book}, or Gromov-Witten invariants.
Their relation to symplectic capacities is not well understood, and we
will not discuss them here.

We start out with a brief description of some relations of symplectic
geometry to neighbouring fields.

\section{Symplectic geometry and its neighbours}  \label{s:neighbours}

Symplectic geometry is a rather new and vigorously developing
mathematical discipline. The ``symplectic explosion`` is described in
\cite{Eliashberg}. Examples of symplectic manifolds are
open subsets of $\left( \R^{2n}, \omega_0 \right)$,
the torus $\R^{2n}/ \Z^{2n}$ endowed with the induced
symplectic form,
surfaces equipped with an area form,
K\"ahler manifolds like complex projective space $\C \P^n$
endowed with their K\"ahler form,
and cotangent bundles with their canonical symplectic
form.
Many more examples are obtained by taking products and through more
elaborate constructions, such as the symplectic blow-up operation.
A diffeomorphism $\varphi$ on a symplectic manifold $(M, \omega)$
is called {\em symplectic}\, or a {\em symplectomorphism}\, if
$\varphi^* \omega = \omega$.

A fascinating feature of symplectic geometry is that it lies at
the crossroad of many other mathematical disciplines.
In this section we mention a few examples of such interactions.

\medskip
{\bf Hamiltonian dynamics.}
Symplectic geometry originated in Hamiltonian dynamics, which
originated in celestial mechanics.
A time-dependent Hamiltonian function on a symplectic manifold $(M,
\omega)$ is a smooth function $H \colon \R \times M \to \R$. Since
$\omega$ is non-degenerate, the equation
\[
\omega \left( X_H, \cdot \right) \,=\, dH ( \cdot )
\]
defines a time-dependent smooth vector field $X_H$ on $M$.
Under suitable assumption on $H$, this vector field generates a
family of diffeomorphisms $\varphi_H^t$ called the {\it Hamiltonian
flow} of $H$.
As is easy to see, each map $\varphi_H^t$ is
symplectic. A {\it Hamiltonian diffeomorphism} $\varphi$ on $M$ is a
diffeomorphism of the form $\varphi_H^1$.

Symplectic geometry is the geometry underlying Hamiltonian
systems. It turns out that this geometric approach to Hamiltonian
systems is very fruitful. Explicit examples are discussed in
\S~\ref{s:caps} below.

\medskip
{\bf Volume geometry.}
A volume form $\Omega$ on a manifold $M$ is a top-dimensional nowhere
vanishing differential form, and a diffeomorphism $\varphi$ of
$M$ is {\it volume preserving} if $\varphi^* \Omega = \Omega$.
Ergodic theory studies the properties of volume preserving
mappings. Its findings apply to symplectic mappings.
Indeed, since a symplectic form $\omega$ is non-degenerate, $\omega^n$ is
a volume form, which is preserved under symplectomorphisms.
In dimension $2$ a symplectic form is just a volume form, so
that a symplectic mapping is just a volume preserving mapping.
In dimensions $2n \ge 4$, however,
symplectic mappings are much more special.
A geometric example for this is Gromov's Nonsqueezing Theorem
stated in \S~\ref{ss:gromov-radius}
and a dynamical example is the (partly solved) Arnol'd
conjecture stating that Hamiltonian diffeomorphisms of closed
symplectic manifolds have at least as many fixed points as
smooth functions have critical points.
For another link between ergodic theory and symplectic geometry
see \cite{Polterovich-loops}.

\medskip
{\bf Contact geometry.}
Contact geometry originated in geometrical optics.
A contact manifold $(P,\alpha)$ is a $(2n-1)$-dimensional
manifold $P$ endowed with a $1$-form $\alpha$ such that $\alpha
\wedge (d \alpha)^{n-1}$ is a volume form on $P$.
The vector field $X$ on $P$ defined by $d \alpha (X, \cdot) =0$
and $\alpha (X)=1$ generates the so-called Reeb flow.
The restriction of a time-independent Hamiltonian system to an
energy surface
can sometimes be realized as the Reeb flow on a contact manifold.
Contact manifolds also arise naturally as boundaries of symplectic manifolds.
One can study a contact manifold $(P,\alpha)$ by symplectic
means by looking at its symplectization
$\left( P \times \R , d (e^t \alpha) \right)$, see e.g.\
\cite{Hofer93a, EGH}.

\medskip
{\bf Algebraic geometry.}
A special class of symplectic manifolds are K\"ahler
manifolds. Such manifolds (and, more generally, complex
manifolds) can be studied by looking at holomorphic curves in
them.
M.~Gromov~\cite{Gr} observed that some of the tools used in the
K\"ahler context can be adapted for the study of symplectic
manifolds. One part of his pioneering work has grown into what is now
called Gromov-Witten theory, see e.g.~\cite{MS2004} for an
introduction.

Many other techniques and constructions from
complex geometry are useful in symplectic geometry.
For example, there is a symplectic version of blowing-up, which
is intimately related to the symplectic packing problem,
see \cite{Mc91,MP94} and \ref{ell:dim4} below.
Another example is Donaldson's construction of symplectic
submanifolds \cite{Donaldson}.
Conversely, symplectic techniques proved useful for
studying problems in algebraic geometry such as Nagata's
conjecture
\cite{B-ample,B-ECM,MP94}
and degenerations of algebraic varieties \cite{B-ICM}.

\medskip
{\bf Riemannian and spectral geometry.}
Recall that the differentiable structure
of a smooth manifold $M$ gives rise to a canonical symplectic form on its
cotangent bundle $T^*M$. Giving a Riemannian metric $g$ on $M$ is
equivalent to prescribing its unit cosphere bundle $S_g^*M \subset
T^*M$, and the restriction of the canonical 1-form from $T^*M$
gives $S^*M$ the structure of a contact manifold.
The Reeb flow on $S_g^*M$ is the geodesic flow (free particle
motion).

In a somewhat different direction, each symplectic form $\omega$  on some
manifold $M$ distinguishes the class of Riemannian metrics which are
of the form $\omega ( J \cdot, \cdot )$ for some almost complex structure $J$.

These (and other) connections between symplectic and Riemannian geometry are
by no means completely explored, and we believe there is still
plenty to be discovered here. Here are some examples of known results
relating Riemannian and symplectic aspects of geometry.

\smallskip
{\it 1. Lagrangian submanifolds.}
A middle-dimensional submanifold $L$ of $(M,\omega)$ is called
{\it Lagrangian} if $\omega$ vanishes on $TL$.

{\it (i) Volume.}
Endow complex projective space $\C \P^n$ with the usual
K\"ahler metric and the usual K\"ahler form.
The volume of submanifolds is taken with respect to this
Riemannian metric.
According to a result of Givental-Kleiner-Oh,
the standard $\R \P^n$ in $\C \P^n$ has minimal volume among all
its Hamiltonian deformations \cite{Oh-1990}.
A partial result for the Clifford torus in $\C \P^n$
can be found in \cite{Goldstein-03}.
The torus $S^1 \times S^1 \subset S^2 \times
S^2$ formed by the equators is also volume minimizing among its
Hamiltonian deformations, \cite{IOS}.
If $L$ is a closed Lagrangian submanifold of
$\left( \R^{2n}, \omega_0 \right)$,
there exists according to \cite{Viterbo00} a constant $C$ depending on
$L$ such that
\begin{equation}  \label{e:Lag:Viterbo}
\vol \left( \varphi_H(L) \right) \,\ge\, C \quad \text{ for all
Hamiltonian deformations of $L$}.
\end{equation}

{\it (ii) Mean curvature.}
The mean curvature form of a Lagrangian submanifold $L$
in a K\"ahler-Einstein manifold can be expressed through symplectic
invariants of $L$, see \cite{CG-04}.

{\it 2. The first eigenvalue of the Laplacian.}
Symplectic methods can be used to estimate the first eigenvalue
of the Laplace operator on functions for certain Riemannian
manifolds \cite{Polterovich-eigenvalue}.

\smallskip
{\it 3. Short billiard trajectories.}
Consider a bounded domain $U \subset \R^n$ with smooth boundary.
There exists a periodic billiard trajectory on $\overline{U}$
of length $l$ with
\begin{equation}  \label{est:billiard}
l^n \,\le\, C_n\, \vol (U)
\end{equation}
where $C_n$ is an explicit constant depending only on $n$, see
\cite{Viterbo00,FGS}.

\section{Examples of symplectic capacities}  \label{s:caps}

In this section we give the formal definition of symplectic
capacities, and discuss a number of examples along with sample
applications.

\subsection{Definition} \label{ss:definition}

Denote by $\Symp^{2n}$ the category of all symplectic manifolds of
dimension $2n$, with symplectic embeddings as morphisms. A {\em
symplectic category}\, is a subcategory $\CC$ of $\Symp^{2n}$ such that
$(M,\om)\in\CC$ implies $(M,\alpha\om)\in\CC$ for all
$\alpha>0$.

{\bf Throughout the paper we will use the symbol $\into$ to denote
symplectic embeddings and $\to$ to denote morphisms in the category $\CC$
(which may be more restrictive).}

Let $B^{2n}(r^2)$ be the open ball of radius $r$ in $\R^{2n}$ and
$Z^{2n}(r^2)=B^2(r^2)\times\R^{2n-2}$ the open cylinder (the reason for
this notation will become apparent below).
Unless stated otherwise, open subsets of $\R^{2n}$ are always equipped with the
canonical symplectic form $\om_0 = \sum_{j=1}^n dy_j \wedge dx_j$.
We will suppress the dimension $2n$ when it is clear from the context
and abbreviate
\[
        B:=B^{2n}(1),\qquad Z:=Z^{2n}(1).
\]

Now let $\CC\subset \Symp^{2n}$ be a symplectic category containing the
ball $B$ 
and the cylinder $Z$. A {\em symplectic capacity}\, on $\CC$
is a covariant functor $c$ from $\CC$ to the category
$([0,\infty],\leq)$ (with $a\leq b$ as morphisms) satisfying
\begin{description}
\item[(Monotonicity)] $c(M,\om)\leq c(M',\om')$ if there exists a
morphism $(M,\om)\to(M',\om')$;
\item[(Conformality)] $c(M,\alpha\om)=\alpha\, c(M,\om)$ for
$\alpha>0$;
\item[(Nontriviality)] $0<c(B)$ 
and $c(Z)<\infty$.
\end{description}
Note that the (Monotonicity) axiom just states the functoriality of
$c$. A symplectic capacity is said to be {\em normalized}\, if
\begin{description}
\item[(Normalization)] $c(B)=1$. 
\end{description}
%
As a frequent example we will use the set $\Op^{2n}$ of open subsets
in $\R^{2n}$. We make it into a symplectic category by identifying
$(U,\alpha^2\om_0)$ with the symplectomorphic manifold $(\alpha
U,\om_0)$ for $U\subset\R^{2n}$ and $\alpha>0$.  We agree that the
morphisms in this category shall be symplectic embeddings induced by
{\em global}\, symplectomorphisms of $\R^{2n}$. With this identification,
the (Conformality) axiom above takes the form
\begin{description}
\item[(Conformality)'] $c(\alpha U)=\alpha^2 c(U)$ for $U\in \Op^{2n}$,
$\alpha>0$.
\end{description}

\subsection{Gromov radius \cite{Gr}} \label{ss:gromov-radius}
In view of Darboux's Theorem one can associate with each
symplectic manifold $(M,\omega)$ the numerical invariant
\[
c_B(M,\om) \,:=\, \sup \left\{ \alpha>0 \mid B^{2n}(\alpha)
\into (M,\om) \right\}
\]
called the {\it Gromov radius}\, of $(M,\omega)$, \cite{Gr}.
It measures the symplectic size of $(M,\omega)$ in a geometric
way, and is reminiscent of the injectivity radius of a
Riemannian manifold. Note that it clearly satisfies the (Monotonicity)
and (Conformality) axioms for a symplectic capacity. It is equally
obvious that $c_B(B)=1$.

If $M$ is $2$-dimensional and connected, then $ \pi c_B (M, \omega) =
\int_M \omega$, i.e.\ $c_B$ is proportional to the volume of $M$, see
\cite{Siburg93}.
The following theorem from Gromov's seminal paper \cite{Gr} implies
that in higher dimensions the Gromov radius is an invariant very
different from the volume.

\smallskip

{\bf Nonsqueezing Theorem (Gromov, 1985).}
{\it The cylinder $Z \in \Symp^{2n}$ satisfies
$c_B (Z) = 1$.
}
\smallskip

In particular, the Gromov radius is a normalized symplectic capacity
on $\Symp^{2n}$. Gromov originally obtained this result by studying
properties of moduli spaces of pseudo-holomorphic curves in symplectic
manifolds.

It is important to realize that the existence of at least one
capacity $c$ with $c(B)=c(Z)$ also {\em implies} the Nonsqueezing Theorem.
We will see below that each of the other important techniques
in symplectic geometry (such as variational methods and the
global theory of generating functions) gave rise to the construction
of such a capacity, and hence an independent proof of this fundamental
result.

It was noted in \cite{EH-90a} that the following result, originally
established by Eliashberg and by Gromov using different methods, is
also an easy consequence of the existence of a symplectic capacity.

\smallskip
{\bf Theorem (Eliashberg, Gromov)}
{\it
The group of symplectomorphisms of a symplectic manifold
$(M,\om)$ is closed for the compact-open $C^0$-topology in the group of
all diffeomorphisms of $M$.
}

\subsection{Symplectic capacities via Hamiltonian systems}

The next four examples of symplectic capacities are constructed
via Hamiltonian systems.
A crucial role in the definition or the construction of these
capacities is played by the action functional of classical
mechanics.
For simplicity, we assume that $(M, \omega) =( \R^{2n}, \omega_0)$.
Given a Hamiltonian function $H \colon S^1 \times \R^{2n} \to
\R$
which is periodic in the time-variable $t \in S^1 = \R /
\Z$ and which generates a global flow $\varphi_H^t$,
the action functional on the loop space $C^\infty(S^1, \R^{2n})$
is defined as
\begin{equation}  \label{e:act}
 \AA_H (\gamma) \,=\, \int_\gamma y\, dx - \int_0^1 H \bigl(
      t,\gamma (t) \bigr) dt .
\end{equation}
Its critical points are exactly the $1$-periodic orbits of $\varphi_H^t$.
Since the action functional is neither bounded from above nor from below,
critical points are saddle points.
In his pioneering work \cite{Rabinowitz78,Rabinowitz79},
P.~Rabinowitz designed special minimax principles adapted to the
hyperbolic structure of the action functional to find such
critical points.
We give a heuristic argument why this works.
Consider the space of loops
\[
E \,=\, H^{1/2}(S^1,\R^{2n}) \,=\,
 \left\{ z \in L^2 \left( S^1; \R^{2n} \right)
\,\Bigg|\, \sum_{k \in \Z}  \left| k \right| \left| z_k \right|^2 < \infty
\right\}
\]
where
$z = \sum_{k \in \Z} e^{2\pi ktJ} z_k$, $z_k \in \R^{2n}$,
is the Fourier series of $z$ and $J$ is the standard complex
structure of $\R^{2n} \cong \C^n$.
The space $E$ is a Hilbert space with inner product
\[
\langle z,w \rangle \,=\, \langle z_0,w_0 \rangle + 2\pi \sum_{k \in
\Z} \left| k \right| \langle z_k,w_k \rangle ,
\]
and there is an orthogonal splitting
$E = E^- \oplus E^0 \oplus E^+$, $z = z^- + z^0 + z^+$,
into the spaces of $z \in E$ having nonzero Fourier coefficients
$z_k \in \R^{2n}$ only for $k<0$, $k=0$, $k>0$.
The action functional $\AA_H \colon C^\infty (S^1,\R^{2n}) \to \R$
extends to $E$ as
\begin{equation}  \label{action:E}
\AA_H(z) \,=\, \left(\tfrac 12 \left\| z^+ \right\|^2 - \tfrac 12 \left\| z^-
\right\|^2\right)  -  \int_0^1 H(t,z(t)) \,dt .
\end{equation}
Notice now the hyperbolic structure of the first term
$\AA_0(x)$,
and that the second term is of lower order.
Some of the critical points $z(t) \equiv const$ of $\AA_0$
should thus persist for $H \neq 0$.

\subsubsection{Ekeland-Hofer capacities \cite{EH-90a,EH-90b}} \label{ss:EH}

The first constructions of symplectic capacities via
Hamiltonian systems were carried out by Ekeland and
Hofer~\cite{EH-90a,EH-90b}.
They considered the space $\FF$ of time-independent
Hamiltonian functions $H \colon \R^{2n} \to [0,\infty)$ satisfying
\begin{itemize}
\item[$\bullet$]
 $H|_U \equiv 0$ for some open subset $U \subset \R^{2n}$, and
\item[$\bullet$]
$H(z) = a |z|^2$ for $|z|$ large, where $a>\pi$, $a\not\in \N\pi$.
\end{itemize}
Given $k\in \N$ and $H\in \FF$, apply equivariant minimax
to define the critical value
\[
c_{H,k} \,:=\, \inf \, \left\{ \sup_{\gamma \in \xi}
 A_H(\gamma) \mid \xi \subset E \text{ \rm is $S^1$-equivariant
and }\ind(\xi) \geq k \right\}
\]
of the action functional~\eqref{action:E},
where $\ind(\xi)$ denotes a suitable Fadell-Rabinowitz
index~\cite{FR:78,EH-90b} of the intersection $\xi \cap S^+$ of $\xi$
with the unit sphere $S^+\subset E^+$. The {\em $k^{th}$ Ekeland-Hofer
capacity} $c_k^{\EH}$ on the symplectic category $\Op^{2n}$ is now
defined as
$$
c_k^{\EH}(U) := \inf \left\{ c_{H,k} \mid H \text{ \rm vanishes on
  some neighborhood of } \ol{U} \right\}
$$
if $U \subset \R^{2n}$ is bounded and as
$$
c_k^{\EH}(U) := \sup \left\{ c_k^{\EH}(V) \mid V \subset U \text{ \rm
  bounded} \right\}
$$
in general. It is immediate from the definition that $c^\EH_1\leq
c^\EH_2\leq c^\EH_3 \leq \dots$ form an increasing sequence.
Their values on the ball and cylinder are
$$
        c^\EH_k(B)=\left[\frac{k+n-1}{n}\right]\pi,\qquad
        c^\EH_k(Z)=k\pi,
$$
where $[x]$ denotes the largest integer $\leq x$.
Hence the existence of $c_1^{\EH}$ gives an independent proof
of Gromov's Nonsqueezing Theorem. Using the capacity $c_n^{\EH}$,
Ekeland and Hofer \cite{EH-90b} also proved the following nonsqueezing
result.

\smallskip
{\bf Theorem (Ekeland-Hofer, 1990)}
{\it
The cube $P=B^2(1) \times \dots \times B^2(1) \subset \C^n$ can be
symplectically embedded into the ball $B^{2n}(r^2)$ if and only if
$r^2 \geq n$.
}

Other illustrations of the use of Ekeland-Hofer capacities
in studying embedding problems for ellipsoids and polydiscs
appear in \S~\ref{s:EP}.

\subsubsection{Hofer-Zehnder capacity \cite{HZ90,HZ}}
\label{ss:HZ}

Given a symplectic manifold $(M,\omega)$ we consider the class
$\cs (M)$ of simple Hamiltonian functions $H \colon M \to [0,\infty)$
characterized by the following properties:
\begin{itemize}
\item[$\bullet$]
$H=0$ near the (possibly empty) boundary of $M$;
\item[$\bullet$]
The critical values of $H$ are $0$ and $\max H$.
\end{itemize}
Such a function is called {\it admissible} if the flow $\varphi_H^t$
of $H$ has no non-constant periodic orbits with period $T \le
1$.

The {\it Hofer-Zehnder capacity} $c_{\HZ}$ on $\Symp^{2n}$ is
defined as
\[
c_{\HZ} (M) \,:=\,
  \sup \left\{ \max H \mid H \in \cs (M) \text{ is admissible} \right\}
\]
It measures the symplectic size of $M$ in a dynamical way. Easily
constructed examples yield the inequality $c_\HZ(B) \ge \pi$. In
\cite{HZ90,HZ}, Hofer and Zehnder applied a minimax technique to the
action functional~\eqref{action:E} to show that $c_\HZ (Z) \le
\pi$, so that
\[
c_\HZ (B) \,=\, c_\HZ (Z) \,=\, \pi,
\]
providing another independent proof of the Nonsqueezing Theorem.
Moreover, for every symplectic manifold $(M,\omega)$ the inequality
$\pi c_B(M) \le c_\HZ(M)$ holds.

The importance of understanding the Hofer-Zehnder capacity comes
from the following result proved in \cite{HZ90,HZ}.

\smallskip
{\bf Theorem (Hofer-Zehnder, 1990)}
{\it
Let $H \colon (M,\omega) \to \R$ be a proper autonomous
Hamiltonian.
If $c_{\HZ}(M) < \infty$, then for almost every $c \in H(M)$ the
energy level $H^{-1}(c)$ carries a periodic orbit.
}
\smallskip

Variants of the Hofer-Zehnder capacity which can be used to
detect periodic orbits in a prescribed homotopy class where
considered in \cite{Lu98,Schwarz00}.

\subsubsection{Displacement energy \cite{Hofer90,LM95}}

Next, let us measure the symplectic size of a subset by looking at
how much energy is needed to displace it from itself. Fix a symplectic
manifold $(M,\om)$. Given a compactly supported Hamiltonian $H \colon
[0,1] \times M \to \R$, set
\[
\| H \| \,:=\, \int_0^1
\left( \sup_{x \in M} H(t,x) - \inf_{x \in M} H(t,x) \right) dt .
\]
The {\it energy}\, of a compactly supported Hamiltonian
diffeomorphism $\varphi$ is
\[
E(\varphi) \,:=\, \inf \left\{ \left\| H \right\| \mid
\varphi = \varphi_H^1 \right\} .
\]
The {\it displacement energy} of a subset $A$ of $M$ is now
defined as
\[
e(A,M) \,:=\, \inf \left\{ E(\varphi) \mid \varphi (A) \cap A =
\emptyset \right\}
\]
if $A$ is compact and as
\[
e(A,M) \,:=\, \sup \left\{ e(K,M) \mid K \subset A \text{ is compact} \right\}
\]
for a general subset $A$ of $M$.

Now consider the special case $(M,\om)=(\R^{2n},\om_0)$. Simple explicit
examples show $e(Z, \R^{2n}) \le \pi$. In \cite{Hofer90}, H.~Hofer
designed a minimax principle for the action functional~\eqref{action:E}
to show that $e(B, \R^{2n}) \ge \pi$, so that
\[
e(B, \R^{2n}) \,=\, e(Z, \R^{2n}) \,=\, \pi .
\]
It follows that $e(\cdot ,\R^{2n})$ is a symplectic capacity on the symplectic
category $\Op^{2n}$ of open subset of $\R^{2n}$.

One important feature of the displacement energy is the inequality
\begin{equation} \label{eq:EC-ineq}
c_{\HZ}(U) \le e(U,M)
\end{equation}
holding for open subsets of many (and possibly all) symplectic
manifolds,
including $(\R^{2n},\om_0)$. Indeed, this inequality and the Hofer-Zehnder Theorem
imply existence of periodic orbits on almost every energy surface of
any Hamiltonian with support in $U$ provided only that $U$ is
displaceable in $M$. The proof of this inequality uses the spectral
capacities introduced in \S~\ref{ss:speccap} below.

As a specific application, consider a closed Lagrangian submanifold
$L$ of $(\R^{2n}, \omega_0 )$. Viterbo~\cite{Viterbo00} used an
elementary geometric construction to show that
\[
e \left( L, \R^{2n} \right) \,\le\, C_n \left( \vol (L)
\right)^{2/n}
\]
for an explicit constant $C_n$. By a result of
Chekanov~\cite{Chekanov1998},
$e \left( L, \R^{2n} \right) >0$.
Since $e \left( \varphi_H(L), \R^{2n} \right) = e \left( L,
\R^{2n} \right)$ for every Hamiltonian diffeomorphism of $L$, we
obtain Viterbo's inequality~\eqref{e:Lag:Viterbo}.

\subsubsection{Spectral capacities \cite{FS,Hofer93,HZ,Oh3,Oh4,Oh5,Sch,V}}
\label{ss:speccap}

For simplicity, we assume again $(M, \omega) =( \R^{2n}, \omega_0)$.
Denote by $\HH$ the space of
compactly supported Hamiltonian functions $H \colon S^1 \times
\R^{2n} \to \R$.
An {\em action selector}\, $\sigma$ selects for each $H \in \HH$
the action $\sigma (H) = \AA_H(\gamma)$ of a ``topologically
visible'' 1-periodic orbit $\gamma$ of $\varphi_H^t$ in a suitable way.
Such {\em action selectors} were constructed by
Viterbo~\cite{V}, who applied minimax to generating functions, and
by Hofer and Zehnder~\cite{Hofer93,HZ}, who applied minimax
directly to the action functional~\eqref{action:E}.
An outline of their constructions can be found in \cite{FGS}.

Given an action selector $\sigma$ for $(\R^{2n},\omega_0)$, one
defines the {\it spectral capacity}\, $c_\sigma$ on the symplectic category
$\Op^{2n}$ by
\[
c_\sigma (U) \,:=\,
\sup \left\{ \sigma(H) \mid H \text{ is supported in } S^1
\times U \right\} .
\]
It follows from the defining properties of an action selector (not
given here) that $c_\HZ(U) \le c_\sigma(U)$ for any spectral capacity
$c_\sigma$. Elementary considerations also imply $c_\sigma(U) \le e
(U,\R^{2n})$, see \cite{FGS,Hofer93,HZ,V}. In this way one in
particular obtains the important inequality (\ref{eq:EC-ineq})
for $M = \R^{2n}$.

Another application of action selectors is

{\bf Theorem (Viterbo, 1992)}
{\it Every non-identical compactly supported Hamiltonian
  diffeomorphism of $\left( \R^{2n}, \omega_0 \right)$ has infinitely
  many non-trivial periodic points.
}

Moreover, the existence of an action selector is an important
ingredient in Viterbo's proof of the estimate~\eqref{est:billiard}
for billiard trajectories.

Using the Floer homology of $(M,\omega)$ filtered by the action
functional, an action selector can be constructed for many (and
conceivably for all) symplectic manifolds $(M,\omega)$,
\cite{FS,Oh3,Oh4,Oh5,Sch}. This existence result implies the
energy-capacity inequality (\ref{eq:EC-ineq}) for arbitrary open
subsets $U$ of such $(M,\om)$, which has many applications~\cite{Sch1}.

\subsection{Lagrangian capacity \cite{CM-03a}}

In~\cite{CM-03a} a capacity is defined on the category of
$2n$-dimensional symplectic manifolds $(M,\om)$ with
$\pi_1(M)=\pi_2(M)=0$ (with symplectic embeddings as morphisms) as
follows. The {\em minimal symplectic area}\,
of a Lagrangian submanifold $L\subset M$ is
$$
        A_{\min}(L) \,:=\, \inf \left\{ \int_\sigma\om \,\Big|\,
        \sigma\in\pi_2(M,L),
        \int_\sigma\om>0 \right\} \,\in\, [0,\infty].
$$
The {\em Lagrangian capacity}\, of $(M,\om)$ is defined as
$$
        c_L(M,\om) \,:=\, \sup \left\{ A_{\min}(L) \mid L\subset M
        \text{ is an embedded Lagrangian torus} \right\}.
$$
Its values on the ball and cylinder are
$$
        c_L(B)=\pi/n,\qquad c_L(Z)=\pi.
$$
As the cube $P=B^2(1) \times \dots \times B^2(1)$ contains the
standard Clifford torus $T^n \subset \C^n$, and is contained in the
cylinder $Z$,
it follows that $c_L(P)=\pi$. Together with
$c_L(B)=\pi/n$ this gives an alternative proof of the nonsqueezing
result of Ekeland and Hofer mentioned in \S~\ref{ss:EH}.
There are also applications of the Lagrangian capacity to Arnold's
chord conjecture and to Lagrangian (non)embedding results into uniruled
symplectic manifolds~\cite{CM-03a}.


\section{General properties and relations between symplectic
capacities}  \label{s:relations}

In this section we study general properties of and relations between
symplectic capacities. We begin by introducing some more notation.
Define the {\em ellipsoids} and {\em polydiscs}
\begin{align*}
        E(a) &:= E(a_1,\dots,a_n) :=
        \left\{ z\in\C^n\;\bigg|\;\frac{|z_1|^2}{a_1}+\dots +
        \frac{|z_n|^2}{a_n}<1 \right\} \cr
        P(a) &:= P(a_1,\dots,a_n) := B^2(a_1)\times\dots\times
        B^2(a_n)
\end{align*}
for $0<a_1\leq\dots\leq a_n\leq\infty$. Note that in
this notation the ball, cube and
cylinder are $B=E(1,\dots,1)$,
$P=P(1,\dots,1)$  and $Z=E(1,\infty,\dots,\infty) = P(1, \infty,
\dots, \infty)$.

Besides $\Symp^{2n}$ and $\Op^{2n}$, two symplectic categories that
will frequently play a role below are
\begin{description}
\item[$\Ell^{2n}$:] the category of ellipsoids in $\R^{2n}$, with
symplectic embeddings induced by global symplectomorphisms of $\R^{2n}$
as morphisms,
\item[$\Pol^{2n}$:] the category of polydiscs in $\R^{2n}$, with
symplectic embeddings induced by global symplectomorphisms of $\R^{2n}$
as morphisms.
\end{description}
\comment{
\begin{remark}
{\rm
The bounded objects $E(a)$ in $\LinEll^{2n}$ and $\Ell^{2n}$ are
the sublevel sets $\left\{ z \in \C^n \mid q_a (z) < 1 \right\}$
of the positive definite quadratic forms
\[
q_a (z) \,=\, \sum_{j=1}^n \frac{|z_j|^2}{a_j} .
\]
This is not a restrictive choice:
Given any positive definite quadratic form $q$ on $\C^n$, there
exists a linear symplectomorphism $A$ of $\C^n$ such that
$q \circ A = q_a$ for some $a \in \R^n$, see~\cite[Section~1.7]{HZ}.
In other words,
\[
\left\{ z \in \C^n \mid q(z) < 1 \right\} = A \left( E(a) \right) .
\]
}
\end{remark}
}

\subsection{Generalized symplectic capacities}  \label{ss:cap}

From the point of view of this work, it is convenient to have a more
flexible notion of symplectic capacities, whose axioms were
originally designed to explicitly exclude such invariants as the volume.
We thus define a {\em generalized symplectic capacity} on a symplectic
category $\CC$ as a covariant functor $c$ from $\CC$ to the category
$([0,\infty],\leq)$ satisfying only the (Monotonicity) and
(Conformality) axioms of \S~\ref{ss:definition}.

Now examples such as the {\em volume capacity} on $\Symp^{2n}$ are
included into the discussion. It is defined as
$$
        c_\vol(M,\om):=\left(\frac{\vol(M,\om)}{\vol(B)}\right)^{1/n},
$$
where $\vol(M,\om):=\int_M\om^n/n!$ is the symplectic volume. For
$n\geq 2$ we have $c_\vol(B)=1$ and $c_\vol(Z)=\infty$, so $c_\vol$ is
a normalized generalized capacity but not a capacity. Many more
examples appear below.

\subsection{Embedding capacities}

Let $\CC$ be a symplectic category. Every object $(X,\Om)$ of $\CC$
induces two generalized symplectic capacities on $\CC$,
\begin{align*}
c_{(X,\Om)}(M,\om) &\,:=\, \sup \left\{ \alpha>0 \mid
(X,\alpha\Om)\to(M,\om) \right\},
\cr
c^{(X,\Om)}(M,\om) &\,:=\, \inf \left\{ \alpha>0 \mid
(M,\om)\to(X,\alpha\Om) \right\},
\end{align*}
Here the supremum and infimum over the empty set are set to $0$ and
$\infty$, respectively. Note that
\begin{equation}\label{eq:inverse}
c_{(X,\Om)}(M,\om) = \left( c^{(M,\om)}(X,\Om) \right)^{-1}.
\end{equation}

\begin{example}  \label{ex:nullinf}
{\rm
Suppose that $(X,\alpha \Om) \to (X,\Om)$ for some $\alpha >1$. Then
$c_{(X,\Om)}(X,\Om)=\infty$ and $c^{(X,\Om)}(X,\Om)=0$, so that
\begin{equation*}
c_{(X,\Om)}(M,\om) \,=\, \left\{
\begin{array}{cl}
\infty & \text{\rm if} \quad (X, \beta \Om) \to (M, \om)
                       \,\text{ \rm for some } \beta >0,\\[.2em]
0 & \text{\rm if} \quad (X, \beta \Om) \to (M, \om)
                       \,\text{ \rm for no } \beta >0,
\end{array} \right.
\end{equation*}
\begin{equation*}
c^{(X,\Om)}(M,\om) \,=\, \left\{
\begin{array}{cl}
0      & \text{\rm if} \quad (M,\om) \to (X, \beta \Om)
                       \,\text{ \rm for some } \beta >0,\\[.2em]
\infty & \text{\rm if} \quad (M,\om) \to (X, \beta \Om)
                       \,\text{ \rm for no } \beta >0.
\end{array} \right.
\end{equation*}
}
\end{example}
The following fact follows directly from the definitions.

\begin{fact}\label{fact:emb1}
Suppose that there exists no morphism $(X,\alpha\Om)\to(X,\Om)$ for
any $\alpha>1$. Then
$c_{(X,\Om)}(X,\Om)=c^{(X,\Om)}(X,\Om)=1$, and for every generalized
capacity $c$ with $0<c(X,\Om)<\infty$,
$$
        c_{(X,\Om)}(M,\om)\leq \frac{c(M,\om)}{c(X,\Om)}\leq
        c^{(X,\Om)}(M,\om)\qquad\text{for all $(M,\om)\in\CC$}.
$$
In other words, $c_{(X,\Om)}$ (resp.\ $c^{(X,\Om)}$) is the minimal
(resp.\ maximal) generalized capacity $c$ with $c(X,\Om) =1$.
\end{fact}

Important examples on $\Symp^{2n}$ arise from the ball $B=B^{2n}(1)$
and cylinder $Z=Z^{2n}(1)$. By Gromov's Nonsqueezing Theorem and
volume reasons we have for $n\geq 2$:
$$
        c_B(Z)=1,\qquad c^Z(B)=1,\qquad c^B(Z)=\infty,\qquad
        c_Z(B)=0.
$$
In particular, for every normalized symplectic capacity $c$,
\begin{equation}  \label{e:ccz}
        c_B(M,\om)\leq c(M,\om)\leq
        c(Z)\,c^Z(M,\om)\qquad\text{for all $(M,\om)\in\Symp^{2n}$}.
\end{equation}
Recall that the capacity $c_B$ is the Gromov radius defined in
\S~\ref{ss:gromov-radius}.
The capacities $c_B$ and $c^Z$ are not comparable on $\Op^{2n}$:
Example~\ref{ex:curious} below shows that for every $k \in \N$ there is a
bounded starshaped domain $U_k$ of $\R^{2n}$ such that
\[
c_B \left( U_k\right) \,\le\, 2^{-k}
\qquad and \qquad
c^Z \left( U_k\right) \,\ge\, \pi k^2,
\]
see also \cite{Hermann}.

We now turn to the question which capacities can be represented as
{\em embedding capacities} $c_{(X,\Om)}$ or $c^{(X,\Om)}$.

\begin{example}
{\rm
Consider the subcategory $\CC \subset \Op^{2n}$ of connected open
sets. Then every generalized capacity $c$ on $\CC$ can be represented
as the capacity $c^{(X,\Om)}$ of embeddings into a (possibly
uncountable) union $(X,\Om)$ of objects in $\CC$.

For this, just define $(X,\Om)$ as the disjoint union of all
$(X_\iota,\Om_\iota)$ in the category $\CC$ with
$c(X_\iota,\Om_\iota)=0$ or $c(X_\iota,\Om_\iota)=1$.
}
\end{example}

\begin{problem}  \label{prob:rep}
Which (generalized) capacities can be represented as $c^{(X,\Om)}$ for a
{\em connected} symplectic manifold $(X,\Om)$?
\end{problem}

\begin{problem}
Which (generalized) capacities can be represented as the capacity
$c_{(X,\Om)}$ of embeddings {\em from} a symplectic manifold $(X,\Om)$?
\end{problem}

\begin{example}  \label{ex:curious}
{\rm
Embedding capacities give rise to some curious generalized
capacities. For example, consider the capacity $c^Y$ of embeddings
into the symplectic manifold $Y:=\amalg_{k\in\N}B^{2n}(k^2)$. It only
takes values $0$ and $\infty$, with $c^Y(M,\om)=0$ iff $(M,\om)$
embeds symplectically into $Y$, cf.\ Example~\ref{ex:nullinf}.
If $M$ is connected,
$\vol(M,\om)=\infty$ implies $c^Y(M,\om)=\infty$. On the other hand,
for every $\eps>0$ there exists an open subset $U\subset\R^{2n}$,
diffeomorphic to a ball, with $\vol(U)<\eps$ and $c^Y(U)=\infty$. To
see this, consider for $k\in\N$ an open neighbourhood $U_k$ of volume
$<2^{-k}\eps$ of the linear cone over the Lagrangian torus $\p
B^2(k^2)\times \dots\times \p B^2(k^2)$. The Lagrangian capacity of $U_k$
clearly satisfies $c_L(U_k)\geq\pi k^2$. The open set $U:=\cup_{k\in\N}U_k$
satisfies $\vol(U)<\eps$ and $c_L(U)=\infty$, hence $U$ does not embed
symplectically into any ball. By appropriate choice of the $U_k$ we
can arrange that $U$ is diffeomorphic to a ball, cf.\
\cite[Proposition~A.3]{Sch}.
}
\diam
\end{example}

\subsubsection*{Special embedding spaces.}  \label{p:special}
Given an arbitrary pair of symplectic manifolds $(X, \Om)$ and $(M,\om)$,
it is a difficult problem to determine or even estimate $c_{(X,\Om)}(M,\om)$
and $c^{(X,\Om)}(M,\om)$.
We thus consider two special cases.

{\bf 1. Embeddings of skinny ellipsoids.}
Assume that $(M,\om)$ is an
ellipsoid $E(a,\dots,a,1)$ with $0< a \leq 1$, and that
$(X,\Om)$ is connected and has finite volume.
Upper bounds for the function
\[
e^{(X,\Om)}(a) \,=\,  c^{(X,\Om)} \left( E(a,\dots,a,1) \right) ,
\quad\, a \in (0,1] ,
\]
are obtained from symplectic embedding results of ellipsoids into $(X,
\Om)$, and lower bounds are obtained from computing other (generalized)
capacities and using Fact~\ref{fact:emb1}.  In particular, the volume
capacity yields
\begin{equation*}
\frac{\left( e^{(X,\Om)}(a) \right)^n}{a^{n-1}} \,\ge\,
\frac{\vol(B)}{\vol (X,\Om)} .
\end{equation*}
The only known general symplectic embedding results for ellipsoids are
obtained  via multiple symplectic folding.
The following result is part of Theorem~3 in \cite{Sch},
which in our setting reads
\begin{fact}  \label{p:limit}
Assume that $(X, \Om)$ is a connected $2n$-dimensional symplectic manifold
of finite volume.
Then
\[
\lim_{a \to 0} \frac{\left( e^{(X,\Om)}(a) \right)^n}{a^{n-1}} \,=\,
\frac{\vol(B)}{\vol (X,\Om)} .
\]
\end{fact}
For a restricted class of symplectic manifolds,
Fact~\ref{p:limit} can be somewhat improved.
The following result is part of Theorem~6.25 of \cite{Sch}.
\begin{fact}  \label{p:limit:fine}
Assume that $X$ is a bounded domain in $\left( \R^{2n}, \om_0 \right)$ with
piecewise smooth boundary or that $(X, \Om)$ is a compact connected
$2n$-dimensional symplectic manifold.
If $n \le 3$, there exists a constant $C>0$ depending only on $(X, \Om)$ such
that
\[
\frac{\left( e^{(X,\Om)}(a) \right)^n}{a^{n-1}} \,\le\, \frac{\vol(B)}{\vol
(X,\Om) \left( 1 - C a^{1/n} \right)} \quad
\text{ for all }\: a < \frac{1}{C^n} .
\]
\end{fact}
These results have their analogues for polydiscs $P(a,\dots,a,1)$.
The analogue of Fact~\ref{p:limit:fine} is known in all dimensions.

{\bf 2. Packing capacities.}
Given an object $(X, \Omega)$ of $\CC$ and $k \in \N$, we denote by
$\coprod_k (X, \Omega)$ the disjoint union of $k$ copies of $(X, \Omega)$
and define
\[
c_{(X,\Om;k)}(M,\om) \,:=\,
\sup \left\{ \alpha>0 \;\Bigg|
\coprod_k (X,\alpha\Om)\into(M,\om) \right\}.
\]
If $\vol (X,\Omega)$ is finite, we see as in Fact~\ref{fact:emb1} that
\begin{equation} \label{ineq:pack}
c_{(X,\Om;k)}(M,\om) \,\le\,
\frac{1}{c_\vol \left(\coprod_k (X,\Om) \right)} c_\vol (M, \om) .
\end{equation}
We say that $(M,\omega)$ admits a {\it full $k$-packing} by $(X,\Om)$ if
equality holds in \eqref{ineq:pack}.
For $k_1, \dots, k_n \in \N$ a full $k_1 \cdots k_n$-packing of
$B^{2n}(1)$ by $E\left( \frac{1}{k_1}, \dots, \frac{1}{k_n} \right)$
is given in \cite{T}.
Full $k$-packings by balls and obstructions to full $k$-packings by
balls are studied in
\cite{Biran97,Biran99,Gr,Kruglikov96,MMT,MP94,Sch,T}.

Assume now that also $\vol (M,\om)$ is finite.
Studying the capacity $c_{(X,\Om;k)}(M,\om)$ is equivalent to studying
the {\em  packing number}
\[
p_{(X,\Om;k)}(M,\om) \,=\, \sup_\alpha \frac{\vol \bigl((\coprod_k \left(X,
\alpha \Om \right) \bigr)}{\vol \left( M,\om\right)}
\]
where the supremum is taken over all $\alpha$ for which $\coprod_k \left(X,
\alpha \Om \right)$ symplectically embeds into $(M,\om)$.
Clearly, $p_{(X,\Om;k)}(M,\om) \le 1$, and equality holds iff equality
holds in \eqref{ineq:pack}.
Results in~\cite{MP94} together with the above-mentioned full
packings of a ball by ellipsoids from \cite{T} imply
\begin{fact}  \label{fact:asympt}
If $X$ is an ellipsoid or a polydisc, then
\[
p_{(X,k)} (M,\om) \to 1 \;\text{ as }\, k \to \infty
\]
for every symplectic manifold $(M,\om)$ of finite volume.
\end{fact}

Note that if the conclusion of Fact~\ref{fact:asympt} holds for $X$ and $Y$,
then it also holds for $X \times Y$.

\begin{problem}
For which bounded convex subsets $X$ of $\R^{2n}$ is the conclusion of
Fact~\ref{fact:asympt} true?
\end{problem}

In~\cite{MP94} and~\cite{Biran97,Biran99}, the packing numbers
$p_{(X,k)}(M)$ are computed for $X=B^4$ and $M=B^4$ or $\C P^2$.
Moreover, the following fact is shown in~\cite{Biran97,Biran99}:
\begin{fact}  \label{fact:ident}
If $X = B^4$, then
for every closed connected symplectic $4$-manifold $(M,\om)$ with
$[\om]\in H^2(M;\Q)$
there exists $k_0(M,\om)$ such that
\[
p_{(X,k)} (M,\om) = 1 \;\text{ for all }\, k \ge k_0 (M,\om) .
\]
\end{fact}
\begin{problem}
For which bounded convex subsets $X$ of $\R^{2n}$ and which connected
symplectic manifolds $(M,\om)$ of finite volume is the conclusion of
Fact~\ref{fact:ident} true?
\end{problem}

\subsection{Operations on capacities}  \label{operations}
We say that a function
$f \colon [0,\infty]^n\to[0,\infty]$ is {\em homogeneous}\, and {\em monotone}\, if
\begin{gather*}
        f(\alpha x_1,\dots,\alpha x_n)=\alpha f(x_1,\dots,x_n)\qquad
        \text{for all }\alpha>0, \cr
        f(x_1,\dots,x_i,\dots,x_n)\leq
        f(x_1,\dots,y_i,\dots,x_n)\qquad \text{for }x_i\leq y_i.
\end{gather*}
If $f$ is homogeneous and monotone and $c_1,\dots,c_n$ are
generalized capacities, then $f(c_1,\dots,c_n)$ is again a generalized
capacity. If in addition $0<f(1,\dots,1)<\infty$ and $c_1,\dots,c_n$
are capacities, then $f(c_1,\dots,c_n)$ is a capacity. Compositions
and pointwise limits of homogeneous monotone functions are again
homogeneous and monotone. Examples include $\max(x_1,\dots,x_n)$,
$\min(x_1,\dots,x_n)$, and the weighted (arithmetic, geometric,
harmonic) means
$$
        \lambda_1x_1+\dots+\lambda_nx_n,\qquad x_1^{\lambda_1}\cdots
        x_n^{\lambda_n},\qquad
        \frac{1}{\frac{\lambda_1}{x_1}+\dots+\frac{\lambda_n}{x_n}}
$$
with $\lambda_1,\dots,\lambda_n\geq 0$,
$\lambda_1+\dots+\lambda_n=1$.

There is also a natural notion of convergence of capacities. We say
that a sequence $c_n$ of generalized capacities on $\CC$ {\em
converges pointwise} to a generalized capacity $c$ if $c_n(M,\om)\to
c(M,\om)$ for every $(M,\om)\in\CC$.

These operations yield lots of
dependencies between capacities, and it is natural to look for
generating systems. In a very general form, this can be formulated as
follows.

\begin{problem}
For a given symplectic category $\CC$, find a minimal
generating system $\GG$ for the (generalized)
symplectic capacities on $\CC$. This means that every (generalized)
symplectic capacity on $\CC$ is the pointwise limit of homogeneous
monotone functions of elements in $\GG$, and no proper subcollection of
$\GG$ has this property.
\end{problem}

This problem is already open for $\Ell^{2n}$ and $\Pol^{2n}$.
One may also ask for generating systems allowing fewer operations,
e.g.~only $\max$ and $\min$, or only positive linear combinations. We
will formulate more specific versions of this problem below. The
following simple fact illustrates the use of operations on
capacities.

\begin{fact} \label{fact:limit}
Let $\CC$ be a symplectic category containing $B$ (resp.~$P$).
Then every generalized capacity $c$ on $\CC$ with $c(B) \neq 0$
(resp.~$c(P) \neq 0$) is the pointwise limit of capacities.
\end{fact}

Indeed, if $c(B)\neq 0$ (resp.~$c(P)\neq 0$), then $c$ is the pointwise
limit as $k\to\infty$ of the capacities
$$
   c_k \,=\, \min \left( c, k\,c_B \right) \Bigl(\text{resp. }
   \min \left( c, k\,c_P \right)  \Bigr).
$$

\begin{example}
{\rm
{\bf (i)}
The generalized capacity $c \equiv 0$ on $\Op^{2n}$ is not a
pointwise limit of capacities, and so the assumption $c(B) \neq 0$
in Fact~\ref{fact:limit} cannot be omitted.

{\bf (ii)}
The assumption $c(B) \neq 0$ is not always necessary:

(a) Define a generalized capacity $c$ on $\Op^{2n}$ by
\begin{equation*}
c (U) \,=\, \left\{
\begin{array}{ll}
0             & \text{\rm if } \vol (U) < \infty ,\\
c_B(U) & \text{\rm if } \vol (U) = \infty.
\end{array} \right.
\end{equation*}
Then $c(B) =0$ and
$c(Z) =1$,
and $c$ is the pointwise limit
of the capacities
$$
c_k \,=\, \max \left( c, \tfrac 1k c_B \right) .
$$
(b) Define a generalized capacity $c$ on $\Op^{2n}$ by
\begin{equation*}
c (U) \,=\, \left\{
\begin{array}{ll}
0       & \text{\rm if } c_B (U) < \infty ,\\
\infty  & \text{\rm if } c_B (U) = \infty.
\end{array} \right.
\end{equation*}
Then $c(B) =0 =c(Z)$ and $c(\R^{2n}) =\infty$, and
$c = \lim_{k \to \infty} \frac 1k c_B$.

{\bf (iii)}
We do not know whether the generalized capacity $c_{\R^{2n}}$ on
$\Op^{2n}$ is the pointwise limit of capacities.
}
\end{example}

\begin{problem}
Given a symplectic category $\CC$ containing $B$ or $P$ and $Z$,
characterize the generalized capacities which are pointwise limits of
capacities.
\end{problem}

\subsection{Continuity}  \label{ss:continuity}

There are several notions of continuity for capacities on open subsets
of $\R^{2n}$, see \cite{Bates95,EH-90a}. For example, consider a {\em
  smooth family of hypersurfaces}
$(S_t)_{-\eps<t<\eps}$ in $\R^{2n}$, each bounding a compact subset with
interior $U_t$. $S_0$ is said to be of {\em restricted contact type}\,
if there exists a vector field $v$ on $\R^{2n}$ which is transverse to
$S_0$ and whose Lie derivative satisfies $L_v\om_0=\om_0$. Let $c$ be
a capacity on $\Op^{2n}$. As the flow of $v$ is conformally
symplectic, the (Conformality) axiom implies (cf.~\cite[p.~116]{HZ})

\begin{fact}\label{fact:cont}
If $S_0$ is of restricted contact type, the function $t\mapsto c(U_t)$
is Lipschitz continuous at $0$.
\end{fact}

Fact~\ref{fact:cont} fails without the hypothesis of restricted
contact type. For example, if $S_0$ possesses no closed characteristic
(such $S_0$ exist by~\cite{G95,G97,GG03}), then by
Theorem~3 in Section~4.2 of~\cite{HZ} the function $t\mapsto
c_\HZ(U_t)$ is not Lipschitz continuous at $0$. V.~Ginzburg~\cite{G03}
presents an example of a smooth family of hypersurfaces
$(S_t)$ (albeit not in $\R^{2n}$) for which the function $t\mapsto
c_\HZ(U_t)$ is not smoother than $1/2$-H\"older continuous.
These considerations lead to

\begin{problem}
Are capacities continuous on all smooth families of domains boun\-ded by
smooth hypersurfaces?
\end{problem}

\subsection{Convex sets}
\label{ss:convex}

Here we restrict to the subcategory $\Conv^{2n} \subset \Op^{2n}$ of
convex open subsets of $\R^{2n}$, with embeddings induced by global
symplectomorphisms of $\R^{2n}$ as morphisms.
Recall that a subset $U\subset \R^{2n}$ is {\it starshaped}\, if
$U$ contains a point $p$ such that for every $q \in U$ the straight
line between $p$ and $q$ belongs to $U$. In particular, convex domains
are starshaped.
\begin{fact}{\rm (Extension after Restriction Principle~\cite{EH-90a})}
\label{f:earp}
Assume that $\varphi \colon U \hookrightarrow \R^{2n}$ is a symplectic
embedding of a bounded starshaped domain $U \subset \R^{2n}$. Then for
any compact subset $K$ of $U$ there exists a symplectomorphism
$\Phi$ of $\R^{2n}$
such that $\Phi |_K = \varphi |_K$.
\end{fact}
This principle continues to hold for some, but not all, symplectic
embeddings of unbounded starshaped domains, see \cite{Sch}.
We say that a capacity $c$ defined on a symplectic
subcategory of $\Op^{2n}$ has the {\em exhaustion property} if
\begin{equation}  \label{def:exh}
c(U) \,=\, \sup \{ \, c(V) \mid V \subset U \text{ \rm is bounded } \}.
\end{equation}
The capacities introduced in \S~\ref{s:caps} all have
this property, but the capacity in Example~\ref{ex:curious} does not.
By Fact~\ref{f:earp}, all statements about capacities defined on a
subcategory of $\Conv^{2n}$ and having the exhaustion property
remain true if we allow all symplectic embeddings (not just those
coming from global symplectomorphisms of $\R^{2n}$) as morphisms.

\begin{fact} \label{fact:emb2}
Let $U$ and $V$ be objects in $\Conv^{2n}$. Then there
exists a morphism $\alpha U \to V$ for every $\alpha\in
(0,1)$ if and only if $c(U)\leq c(V)$ for all generalized
capacities $c$ on $\Conv^{2n}$.
\end{fact}
Indeed, the necessity of the condition is obvious, and the sufficiency
follows by observing that $\alpha U \to U$ for all $\alpha\in (0,1)$
and $1 \leq c_U(U) \leq c_U(V)$. What happens for $\alpha=1$ is not
well understood, see \S~\ref{ss:recog} for related discussions. The
next example illustrates that the conclusion of Fact~\ref{fact:emb2}
is wrong without the convexity assumption.
\begin{example}
{\rm
Consider the open annulus $A = B(4) \setminus B(1)$ in
$\R^2$.
If $\frac 3 4 < \alpha^2< 1 $, then $\alpha A$ cannot be
embedded into $A$ by a global symplectomorphism.
Indeed, volume considerations show that any potential such
global symplectomorphism
would have to map $A$ homotopically nontrivially into itself. This
would force the image of the ball $\alpha B(1)$ to cover all
of $B(1)$, which is impossible for volume reasons.
\diam
}
\end{example}

Assume now that $c$ is a normalized symplectic capacity on $\Conv^{2n}$.
Using John's ellipsoid, Viterbo \cite{Viterbo00} noticed that
there is a constant $C_n$ depending only on $n$ such that
\begin{equation*}
c^Z (U) \,\le\, C_n\, c_B (U) \quad\,\text{for all }\, U \in \Conv^{2n}
\end{equation*}
and so, in view of \eqref{e:ccz},
\begin{equation}
\label{viterbo_estimate}
c_B(U) \,\le\, c(U) \,\le\,
C_n \,c(Z) \,c_B (U) \quad\,\text{for all }\, U \in \Conv^{2n} .
\end{equation}
In fact, $C_n \le (2n)^2$ and $C_n \le 2n$ on centrally symmetric convex sets.

\begin{problem}  \label{problem:Cn}
What is the optimal value of the constant $C_n$ appearing in
(\ref{viterbo_estimate})?
In particular, is $C_n=1$?
\end{problem}
%

Note that $C_n=1$ would imply uniqueness of capacities satisfying
$c(B)=c(Z)=1$ on $\Conv^{2n}$.
In view of Gromov's Nonsqueezing Theorem, $C_n=1$ on $\Ell^{2n}$ and
$\Pol^{2n}$. More generally, this equality holds for all convex
Reinhardt domains \cite{Hermann}. In particular, for these special
classes of convex sets
\[
\pi c_B \,=\, c_1^\EH \,=\, c_\HZ \,=\, e(\cdot ,\R^{2n}) \,=\,  \pi c^Z.
\]

\subsection{Recognition}  \label{ss:recog}

One may ask how complete the information provided by all symplectic
capacities is. Consider two objects $(M,\om)$ and $(X,\Om)$ of a
symplectic category $\CC$.

\begin{question}  \label{q:1}
Assume $c(M,\om) \le c(X,\Om)$ for all generalized symplectic capacities
$c$ on $\CC$. Does it follow that $(M,\om) \into (X,\Om)$ or even that
$(M,\om)  \to (X,\Om)$?
\end{question}

\begin{question}  \label{q:2}
Assume $c(M,\om) = c(X,\Om)$ for all generalized symplectic capacities
$c$ on $\CC$.
Does it follow that $(M,\om)$ is symplectomorphic to $(X,\Om)$ or even
that $(M,\om) \cong (X,\Om)$ in the category $\CC$?
\end{question}

Note that if $(M,\alpha\om) \to (M,\om)$ for all $\alpha\in (0,1)$
then, under the assumptions of Question~\ref{q:1}, the argument
leading to Fact~\ref{fact:emb2} yields $(M,\alpha \om) \to (X,\Om)$
for all $\alpha \in (0,1)$.

\begin{example}  \label{example:recon}
{\rm
{\bf (i)}
Set $U = B^2(1)$ and $V = B^2(1) \setminus \{0\}$.
For each $\alpha <1$ there exists a symplectomorphism of $\R^2$
with $\varphi \left( \alpha U \right) \subset V$, so that
monotonicity and conformality imply $c(U) = c(V)$ for all
generalized capacities $c$ on $\Op^2$.
Clearly, $U \hookrightarrow V$, but
$U \nrightarrow V$, and $U$ and $V$ are not symplectomorphic.

{\bf (ii)}
Set $U = B^2(1)$ and let
$V = B^2(1) \setminus \{ (x,y) \mid x \ge 0,\, y=0 \}$ be the
slit disc.
As is well-known, $U$ and $V$ are symplectomorphic.
Fact~\ref{f:earp} implies $c(U) = c(V)$ for all generalized
capacities $c$ on $\Op^2$, but clearly $U \nrightarrow
V$.
In dimensions $2n \ge 4$ there are bounded convex sets $U$ and
$V$ with smooth boundary which are symplectomorphic while $U
\nrightarrow V$, see \cite{ElH}.

{\bf (iii)}
Let $U$ and $V$ be ellipsoids in $\Ell^{2n}$.
The answer to Question~\ref{q:1} is unknown even for $\Ell^4$.
For $U = E(1,4)$ and $V=B^4(2)$ we have $c(U) \le c(V)$ for all
generalized capacities that can presently be
computed, but it is unknown whether $U \hookrightarrow V$,
cf.~\ref{ell:dim4} below.
By Fact~\ref{fact:EHdetermine} below, the answer to Question~\ref{q:2}
is ``yes'' on $\Ell^{2n}$.

{\bf (iv)}
Let $U$ and $V$ be polydiscs in $\Pol^{2n}$.
Again, the answer to Question~\ref{q:1} is unknown even for
$\Pol^4$. However, in this dimension the Gromov radius together with
the volume capacity determine a polydisc, so that the answer to
Question~\ref{q:2} is ``yes'' on $\Pol^4$.
\diam
}
\end{example}

\smallskip
\begin{problem}
Are two polydiscs in dimension $2n \ge 6$ with equal
generalized symplectic capacities symplectomorphic?
\end{problem}

To conclude this section, we mention a specific example in which $c(U)
= c(V)$ for all known (but possibly not for all) generalized
symplectic capacities.

\begin{example}
{\rm
Consider the subsets
\[
U = E(2,6) \times E(3,3,6) \quad \text{ and }
V= E(2,6,6) \times E(3,3)
\]
of $\R^{10}$.
Then $c(U) =c(V)$ whenever $c(B)=c(Z)$ by the Nonsqueezing
Theorem, the volumina agree, and $c_k^\EH (U) = c_k^{EH}(V)$ for
all $k$ by the product formula (\ref{e:Chekanov}).
It is unknown whether $U \hookrightarrow V$ or $V
\hookrightarrow U$ or $U \to V$.
Symplectic homology as constructed in \cite{FH94,Traynor94}
does not help in these problems because a computation based on
\cite{FHW94} shows that all symplectic homologies of $U$ and $V$
agree.
}
\end{example}

\subsection{Hamiltonian representability}  \label{ss:repres}

Consider a bounded domain $U \subset \R^{2n}$ with smooth
boundary of restricted contact type (cf. \S~\ref{ss:continuity} for
the definition).
A {\it closed characteristic} $\gamma$ on $\partial U$
is an embedded circle in $\partial U$ tangent to the
characteristic line bundle
\[
\LL_U \,=\, \left\{ (x, \xi) \in T \:\! \partial U \mid \om_0(\xi, \eta) =0
\;\text{ for all }\, \eta \in T_x \:\! \partial U \right\} .
\]
If $\partial U$ is represented as a regular energy surface
$\left\{ x \in \R^{2n} \mid H(x) = \const \right\}$
of a smooth function $H$ on $\R^{2n}$,
then the Hamiltonian vector field $X_H$ restricted to $\partial U$ is
a section of $\LL_U$, and so the traces of the periodic
orbits of $X_H$ on $\partial U$ are the closed characteristics on
$\partial U$.
\index{action of a closed characteristic}
The {\it action} \:\!$A \left( \gamma \right)$
of a closed characteristic $\gamma$ on $\partial U$ is defined as
$A \left( \gamma \right) = \left| \int_\gamma y\,dx \right|$.
The set
\[
\Sigma \left( U \right) \,=\,
\bigl\{ k A \left( \gamma \right) \mid k = 1,2, \dots ; \;
\gamma \text{ is a closed characteristic on } \partial U \bigr\}
\]
is called the {\it action spectrum} of $U$.
This set is nowhere dense in $\R$, cf.~\cite[Section~5.2]{HZ},
and it is easy to see that $\Sigma (U)$ is closed and $0 \notin
\Sigma (U)$.
For many capacities $c$ constructed via Hamiltonian systems, such as
Ekeland-Hofer capacities $c_k^\EH$ and spectral capacities $c_\sigma$,
one has $c(U) \in \Sigma (U)$, see \cite{EH-90b,He-inner}.
Moreover,
\begin{equation}  \label{id:convex}
c_\HZ (U) \,=\, c^\EH_1(U) \,=\, \min \left( \Sigma (U)\right)
\quad \text{ if $U$ is convex.}
\end{equation}
One might therefore be tempted to ask

\begin{question}  \label{q:hermann}
Is it true that $\pi c(U) \in \Sigma (U)$ for every normalized
symplectic capacity $c$ on $\Op^{2n}$ and every domain $U$ with
boundary of restricted contact type?
\end{question}

The following example due to D.~Hermann~\cite{Hermann} shows that
the answer to Question~\ref{q:hermann} is ``no''.

\begin{example}\label{ex:hermann}
{\rm
Choose any $U$ with boundary of restricted contact type such
that
\begin{equation}  \label{ass}
c_B (U) \,<\, c^Z(U) .
\end{equation}
Examples are bounded starshaped domains $U$ with smooth boundary
which contain the Lagrangian torus $S^1 \times \dots \times S^1$
but have small volume:
According to \cite{Sikorav91}, $c^Z (U) \ge 1$, while $c_B(U)$
is as small as we like.
Now notice that for each $t \in [0,1]$,
\[
c_t \,=\, (1-t) c_B + t c^Z
\]
is a normalized symplectic capacity on $\Op^{2n}$. By
\eqref{ass}, the interval
\[
\left\{ c_t (U) \mid t \in [0,1] \right\}
\,=\, [ c_B(U), c^Z (U)]
\]
has positive measure and hence cannot lie
in the nowhere dense set $\Sigma (U)$.
\diam
}
\end{example}

D.~Hermann also pointed out that the argument in
Example~\ref{ex:hermann} together with \eqref{id:convex} implies
that the question ``$C_n=1$?'' posed in Problem~\ref{problem:Cn} is
equivalent to Question~\ref{q:hermann} for convex sets.

\subsection{Products}  \label{products}

Consider a family of symplectic categories $\CC^{2n}$ in all
dimensions $2n$ such that
$$
   (M,\om)\in\CC^{2m},\ \ (N,\sigma)\in\CC^{2n} \;\Longrightarrow\;
   (M\times N,\om\oplus\sigma)\in\CC^{2(m+n)}.
$$
We say that a collection
$c \colon \amalg_{n=1}^\infty\CC^{2n}\to[0,\infty]$ of generalized capacities
has the {\em product property}\, if
$$
   c(M\times N,\om\oplus\sigma) \,=\, \min \{c(M,\om),c(N,\sigma)\}
$$
for all $(M,\om)\in\CC^{2m}$, $(N,\sigma)\in\CC^{2n}$. If
$\R^2 \in \CC^2$ and $c(\R^2)=\infty$, the product property implies the
{\em stability property}
$$
   c(M\times\R^2,\om\oplus\om_0) \,=\, c(M,\om)
$$
for all $(M,\om)\in\CC^{2m}$.

\begin{example}  \label{ex:product}
{\rm
{\bf (i)}
Let $\Sigma_g$ be a closed surface of genus $g$ endowed with an area
form $\omega$. Then
\begin{equation*}
c_B \left( \Sigma_g \times \R^2, \omega \oplus \omega_0 \right) \,=\,
\left\{
\begin{array}{lll}
c_B \left( \Sigma_g,\omega \right) = \frac{1}{\pi} \omega \left( \Sigma_g \right)
      &  \text{\rm if} & g=0 , \\[.5em]
   \infty & \text{\rm if} & g \ge 1 .
\end{array} \right.
\end{equation*}
While the result for $g=0$ follows from Gromov's Nonsqueezing Theorem,
the result for $g \ge 1$ belongs to
Polterovich~\cite[Exercise~12.4]{MS} and Jiang~\cite{Jiang00}.
Since $c_B$ is the smallest normalized symplectic capacity on
$\Symp^{2n}$, we find that no collection $c$ of symplectic capacities
defined on the family $\coprod_{n=1}^\infty \Symp^{2n}$
with $c \left( \Sigma_g, \omega \right) < \infty$ for some $g
\ge 1$ has the product or stability property.

{\bf (ii)}
On the family of polydiscs $\coprod_{n=1}^\infty \Pol^{2n}$,
the Gromov radius, the Lagrangian capacity and the unnormalized
Ekeland-Hofer capacities $c_k^{\EH}$ all have the product property
(see Section~\ref{ss:poly}).
The volume capacity is not stable.

{\bf (iii)}
Let $U \in \Op^{2m}$ and $V \in \Op^{2n}$ have smooth boundary of
restricted contact type (cf.~\S~\ref{ss:continuity} for the
definition). The formula
\begin{equation}  \label{e:Chekanov}
c_k^\EH \left( U \times V \right) \,=\, \min_{i+j=k} \left( c_i^\EH
\left(U\right) + c_j^\EH \left(V\right) \right),
\end{equation}
in which we set $c^\EH_0 \equiv 0$,
was conjectured by Floer and Hofer \cite{Viterbo89}
and has been proved by Chekanov \cite{Chekanov2008} as an application
of his equivariant Floer homology.
Consider the collection of sets $U_1 \times \dots \times U_l$, where
each $U_i \in \Op^{2n_i}$ has smooth boundary of restricted contact
type, and $\sum_{i=1}^l n_i = n$.
We denote by $\RCT^{2n}$ the corresponding category with symplectic
embeddings induced by global symplectomorphisms of $\R^{2n}$ as morphisms.
If $v_i$ are vector fields on $\R^{2n_i}$ with $L_{v_i}\om_0=\om_0$,
then $L_{v_1+ \dots + v_l}\om_0=\om_0$ on $\R^{2n}$. Elements of
$\RCT^{2n}$ can therefore be exhausted by elements of $\RCT^{2n}$ with
smooth boundary of restricted contact type.
This and the exhaustion property~\eqref{def:exh} of the
$c^\EH_k$ shows that \eqref{e:Chekanov} holds for all $U \in \RCT^{2m}$
and $V \in \RCT^{2n}$, implying in particular that Ekeland-Hofer
capacities are stable on $\RCT := \coprod_{n=1}^\infty \RCT^{2n}$. Moreover,
\eqref{e:Chekanov} yields that
\[
c^\EH_k \left( U \times V \right) \,\le\,
\min \left( c_k^\EH \left(U\right),\, c_k^\EH \left(V\right) \right),
\]
and it shows that $c^\EH_1$ on $\RCT$ has the product property.
Using~\eqref{e:Chekanov} together with an induction over the
number of factors and
$c^\EH_2 \left( E(a_1, \dots, a_n) \right) \le 2 a_1$
we also see that $c^\EH_2$ has the product property on products
of ellipsoids.
For $k \ge 3$, however, the Ekeland-Hofer capacities $c^\EH_k$ on $\RCT$
do not have the product property. As an example,
for $U = B^4(4)$ and $V=E(3,8)$ we have
\[
c^\EH_3 (U\times V) = 7 < 8=\min \left( c^\EH_3 (U),c^\EH_3 (V)\right).
\]
}
\end{example}

\begin{problem}
Characterize the collections of (generalized) capacities on polydiscs
that have the product (resp.~stability) property.
\end{problem}

Next consider a collection $c$ of generalized capacities
on open subsets $\Op^{2n}$. In general, it will not be stable.
However, we can stabilize $c$ to obtain stable generalized capacities
$c^\pm \colon \coprod_{n=1}^\infty {\Op}^{2n} \to [ 0, \infty ]$,
\[
   c^+(U) := \limsup_{k\to\infty}c(U\times\R^{2k}),\qquad
   c^-(U) := \liminf_{k\to\infty}c(U\times\R^{2k}).
\]
Notice that $c(U) = c^+(U)=c^-(U)$ for all $U \in \coprod_{n=1}^\infty
\Op^{2n}$ if and only if $c$ is stable.
If $c$ consists of capacities and there exist constants $a,A>0$ such
that
\[
   a\leq c\Bigl(B^{2n}(1)\Bigr)\leq c\Bigl( Z^{2n}(1)\Bigr)\leq
   A\qquad \text{for all }n\in\N,
\]
then $c^\pm$ are collections of capacities. Thus there exist plenty of
stable capacities on $\Op^{2n}$. However, we have
\begin{problem}
Decide stability of specific collections of capacities on $\Conv^{2n}$
or $\Op^{2n}$, e.g.:
Gromov radius, Ekeland-Hofer capacity, Lagrangian capacity,
and the embedding capacity $c_P$ of the unit cube.
\end{problem}
\begin{problem}
Does there exist a collection of capacities on
$\coprod_{n=1}^\infty \Conv^{2n}$ or $\coprod_{n=1}^\infty \Op^{2n}$
with the product property?
\end{problem}

\subsection{Higher order capacities\,?}  \label{higher}
Following \cite{Hofer90b}, we briefly discuss the concept of
higher order capacities. Consider a symplectic category $\CC \subset
\Symp^{2n}$ containing $\Ell^{2n}$ and fix $d \in \{ 1, \dots, n \}$.
A {\em symplectic $d$-capacity}\, on $\CC$
is a generalized capacity satisfying
\begin{description}
\item[($d$-Nontriviality)] $0<c(B)$ and
\begin{equation*}
\left\{
 \begin{array}{ll}
  c \left( B^{2d} (1) \times \R^{2(n-d)} \right) < \infty,\\[.3em]
  c \left( B^{2(d-1)} (1) \times \R^{2(n-d+1)} \right) = \infty.
 \end{array}
\right.
\end{equation*}
\end{description}
For $d=1$ we recover the definition of a symplectic capacity, and for
$d=n$ the volume capacity $c_\vol$ is a symplectic $n$-capacity.
%

\begin{problem}  \label{prob:dcap}
Does there exist a symplectic $d$-capacity on a symplectic category
$\CC$ containing $\Ell^{2n}$ for some $d \in \{2, \dots, n-1 \}$?
\end{problem}
Problem~\ref{prob:dcap} on $\Symp^{2n}$ is equivalent to the
following symplectic embedding problem.
\begin{problem}  \label{prob:emb}
Does there exist a symplectic embedding
\begin{equation}  \label{eq:emb}
B^{2(d-1)} (1) \times \R^{2(n-d+1)} \,\into\, B^{2d}(R) \times \R^{2(n-d)}
\end{equation}
for some $R<\infty$ and $d \in \{2, \dots, n-1 \}$?
\end{problem}
Indeed, the existence of such an embedding would imply that no
symplectic \text{$d$-capacity} can exist on $\Symp^{2n}$. Conversely, if no
such embedding exists, then the embedding capacity $c^{Z_{2d}}$ into
$Z_{2d} = B^{2d} (1) \times \R^{2(n-d)}$ would be an example of
a $d$-capacity on $\Symp^{2n}$.
The Ekeland-Hofer capacity $c_d^{\EH}$ shows that $R \ge 2$ if a
symplectic embedding \eqref{eq:emb} exists.
The known symplectic embedding techniques are not designed to
effectively use the unbounded factor of the target space in
\eqref{eq:emb}. E.g., multiple symplectic folding only shows that
there exists a function $f \colon [1, \infty) \to \R$ with $f(a) <
\sqrt{2a}+2$ such that for each $a \ge 1$ there exists a symplectic
embedding
\[
B^2(1) \times B^2(a) \times \R^2 \,\hookrightarrow\, B^4 \left( f(a)
\right) \times \R^2
\]
of the form $\varphi \times id_2$, see \cite[Section 4.3.2]{Sch}.

\section{Ellipsoids and polydiscs}  \label{s:EP}

In this section we investigate generalized capacities on the
categories of ellipsoids $\Ell^{2n}$ and polydiscs $\Pol^{2n}$ in more
detail.
All (generalized) capacities $c$ in this section are defined on some
symplectic subcategory of $\Op^{2n}$ containing at least one of the above
categories and are assumed to have the exhaustion property~\eqref{def:exh}.

\subsection{Ellipsoids}

\subsubsection{Arbitrary dimension}  \label{ell:arbitrary}

We first describe the values of the capacities introduced in
\S~\ref{s:caps} on ellipsoids.

The values of the Gromov radius $c_B$ on ellipsoids are
$$
        c_B\bigl(E(a_1,\dots,a_n)\bigr)=\min\{a_1,\dots,a_n\}.
$$
More generally, monotonicity implies that this formula holds for
all symplectic capacities $c$ on $\Op^{2n}$ with $c(B) = c(Z)
=1$ and hence also for $\frac{1}{\pi}c_1^\EH$,
$\frac{1}{\pi} c_\HZ$, $\frac{1}{\pi}e(\cdot ,\R^{2n})$ and $c^Z$.

The values of the Ekeland-Hofer capacities on the ellipsoid
$E(a_1,\dots,a_n)$ can be described as follows~\cite{EH-90b}. Write the
numbers $m\,a_i\:\!\pi$, $m\in\N$, $1\leq i\leq n$, in increasing order as
$d_1\leq d_2\leq\dots$, with repetitions if a number occurs several
times. Then
$$
        c^\EH_k\bigl(E(a_1,\dots,a_n)\bigr) = d_k.
$$

The values of the Lagrangian capacity on ellipsoids are presently not
known.
In \cite{CM-03b}, Cieliebak and Mohnke expect to prove the following
\begin{conjecture} \label{conj:c_L}
$$
        c_L\bigl(E(a_1,\dots,a_n)\bigr) =
        \frac{\pi}{1/a_1+\dots+1/a_n}.
$$
\end{conjecture}

Since $\vol\bigl(E(a_1,\dots,a_n)\bigr)=a_1\cdots a_n\vol(B)$, the
values of the volume capacity on ellipsoids are
$$
        c_\vol\bigl(E(a_1,\dots,a_n)\bigr)=(a_1\cdots a_n)^{1/n}.
$$


In view of conformality and the exhaustion property, a
(generalized) capacity on $\Ell^{2n}$ is
determined by its values on
the ellipsoids $E(a_1,\dots,a_n)$ with $0<a_1\leq\dots\leq a_n=1$.
So we can view each (generalized) capacity $c$ on ellipsoids as a function
$$
        c(a_1,\dots,a_{n-1}) := c \left( E(a_1,\dots,a_{n-1},1) \right)
$$
on the set $\{0< a_1\leq\dots\leq a_{n-1}\leq 1\}$.
By Fact~\ref{fact:cont}, this function is continuous. This identification
with functions yields a notion of {\em uniform convergence} for
capacities on $\Ell^{2n}$.

For what follows, it is useful to have normalized versions of the
Ekeland-Hofer capacities, so in dimension $2n$ we define
$$
        \bar c_k:= \frac{c^\EH_k}{[\frac{k+n-1}{n}]\pi}.
$$

\begin{prop}\label{prop:lim_ell}
As $k\to\infty$, for every $n \ge 2$ the normalized Ekeland-Hofer
capacities $\bar c_k$ converge uniformly on $\Ell^{2n}$ to the
normalized symplectic capacity $c_\infty$ given by
$$
c_\infty\left(E(a_1,\dots,a_n)\right) = \frac{n}{1/a_1+\dots+1/a_n}.
$$
\end{prop}

\begin{remark}
Note that Conjecture~\ref{conj:c_L} asserts that $c_\infty$ agrees
with the normalized Lagrangian capacity $\bar c_L= n c_L/\pi$ on
$\Ell^{2n}$.
\end{remark}

{\it Proof of Proposition~\ref{prop:lim_ell}.}
Fix $\eps >0$. We need to show that $\left| \bar c_k (a) - c_\infty (a)
\right| \le \eps$ for every vector $a=(a_1,\dots,a_n)$
with $0<a_1\le a_2\le \dots \le a_n =1$ and all sufficiently large $k$.
Abbreviate $\delta = \eps /n$.

\smallskip
{\em Case 1.}
$a_1 \le \delta$. Then
$$
        c_k^\EH(a)\leq k\delta\pi,\qquad \bar c_k(a)\leq
        n\delta,\qquad c_\infty(a)\leq n\delta
$$
from which we conclude $\left| \bar c_k(a)- c_\infty(a) \right| \le
n\delta = \eps$ for all $k\geq 1$.

\smallskip
{\em Case 2.}
$a_1 > \delta$. Let $k \ge 2 \frac{n-1}{\delta} +2$.
For the unique integer $l$ with
$$
        \pi l\,a_n\leq c^\EH_k(a)<\pi(l+1)a_n
$$
we then have $l \ge 2$.
In the increasing sequence of the numbers $m\,a_i$ ($m\in\N$, $1 \le
i\leq n$), the first $[l\,a_n / a_i]$ multiples of $a_i$ occur no later
than $l\,a_n$. By the description of the Ekeland-Hofer capacities on
ellipsoids given above, this yields the estimates
$$
        \frac{(l-1)\,a_n}{a_1}+\dots+\frac{(l-1)\,a_n}{a_n} \leq k \le
        \frac{(l+1)\,a_n}{a_1}+\dots+\frac{(l+1)\,a_n}{a_n}.
$$
With $\gamma:=a_n/a_1+\dots+a_n/a_n$ this becomes
$$
        (l-1)\gamma\leq k\leq (l+1)\gamma.
$$
Using $\gamma\geq n$, we derive the inequalities
\begin{align*}
        \left[\frac{k+n-1}{n}\right] &\leq \frac{k}{n}+1\leq
        \frac{(l+1)\gamma+n}{n} \leq \frac{(l+2)\gamma}{n}, \cr
        \left[\frac{k+n-1}{n}\right] &\geq \frac{k}{n}\geq
        \frac{(l-1)\gamma}{n}.
\end{align*}
With the definition of $\bar c_k$ and the estimate above for
$c^\EH_k$, we find
$$
        \frac{n\,l\,a_n}{(l+2)\gamma}\leq \bar c_k(a) =
        \frac{c^\EH_k(a)}{[\frac{k+n-1}{n}]\pi} \leq
        \frac{n(l+1)a_n}{(l-1)\gamma}.
$$
Since $c_\infty(a)=n\,a_n/\gamma$, this becomes
$$
        \frac{l}{l+2} c_\infty(a) \leq \bar c_k(a)
        \leq \frac{l+1}{l-1} c_\infty(a),
$$
which in turn implies
$$
        |\bar c_k(a)- c_\infty(a)| \leq \frac{2 c_\infty(a)}{l-1}.
$$
Since $a_1 > \delta$ we have
$$
        \gamma \le \frac{n}{\delta},\qquad l+1 \ge \frac{k}{\gamma} \ge
\frac{k\delta}{n},
$$
from which we conclude
$$
        |\bar c_k(a)- c_\infty(a)| \leq \frac{2}{l-1}\le
        \frac{2n}{k\delta-2n} \le \eps
$$
for $k$ sufficiently large.
\proofend
We turn to the question whether Ekeland-Hofer
capacities generate the space of all capacities on
ellipsoids by suitable operations. First note some easy facts.

\begin{fact}  \label{fact:EHdetermine}
An ellipsoid $E\subset\R^{2n}$ is uniquely determined by its
Ekeland-Hofer capacities $c^\EH_1(E), c^\EH_2(E),\dots$.
\end{fact}

Indeed, if $E(a)$ and $E(b)$ are two ellipsoids with $a_i=b_i$ for
$i<k$ and $a_k<b_k$, then the multiplicity of $a_k$ in the sequence of
Ekeland-Hofer capacities is one higher for $E(a)$ than for $E(b)$, so
not all Ekeland-Hofer capacities agree.

\begin{fact}
For every $k\in\N$ there exist ellipsoids $E$ and $E'$ with
$c^\EH_i(E)=c^\EH_i(E')$ for $i<k$ and $c^\EH_k(E)\neq c^\EH_k(E')$.
\end{fact}

For example, we can take $E=E(a)$ and $E'=E(b)$ with $a_1=b_1=1$,
$a_2=k-1/2$,
$b_2=k+1/2$, and $a_i=b_i=2k$ for $i\geq 3$.
So formally, every generalized capacity on ellipsoids is a function of
the Ekeland-Hofer capacities, and the Ekeland-Hofer capacities are
functionally independent.
However, Ekeland-Hofer capacities do not form a
generating system for symplectic capacities on $\Ell^{2n}$ (see
Example~\ref{ex:333} below), and on bounded ellipsoids
each finite set of Ekeland-Hofer capacities
is determined by the (infinitely many) other Ekeland-Hofer
capacities:
\begin{lemma}  \label{lem:reconstruction}
Let $d_1 \le d_2 \le \dots$ be an increasing sequence of real
numbers obtained from the sequence $c^{\EH}_1 (E) \leq c^{\EH}_2(E)
\leq \dots$ of Ekeland-Hofer capacities of a bounded ellipsoid
$E \in \Ell^{2n}$ by removing at most $N_0$ numbers. Then $E$ can be
recovered uniquely.
\end{lemma}

\begin{proof}
We first consider the special case in which $E=E(a_1,\dots,a_n)$
is such that $a_i/a_j \in \Q$ for all $i,j$.
In this case, the sequence $d_1 \leq
d_2 \leq \dots$ contains infinitely many blocks of $n$ consecutive
equal numbers. We traverse the sequence until we have found $N_0+1$
such blocks, for each block $d_k=d_{k+1}=\dots=d_{k+n-1}$ recording the
number $g_k:=d_{k+n}-d_{k}$. The minimum of the $g_k$ for the
$N_0+1$ first blocks equals $a_1$. After deleting each occurring positive
integer multiple of $a_1$ once from the sequence $d_1\leq d_2 \leq
\dots$, we can repeat the same procedure to determine $a_2$, and so
on.

In general, we do not know whether or not $a_i/a_j \in \Q$ for all $i,j$.
To reduce to the previous case, we split the
sequence $d_1 \leq d_2 \leq \dots$ into (at most $n$)
subsequences of numbers with rational quotients.
More precisely we traverse the sequence, grouping the $d_i$
into increasing subsequences $s_1, s_2,\dots$, where each new number
is added to the first subsequence $s_j$ whose members are
rational multiples of it.
Furthermore, in this process we record for each sequence $s_j$ the
maximal length $l_j$ of a block of consecutive equal numbers seen so far.
We stop when
\begin{itemize}
\item[(i)] the sum of the $l_j$ equals $n$, and
\item[(ii)] each subsequence $s_j$ contains at least $N_0+1$ blocks of
$l_j$ consecutive equal numbers.
\end{itemize}
Now the previously described procedure in the case that $a_i/a_j
\in \Q$ for all $i,j$
can be applied for each subsequence $s_j$ separately, where $l_j$
replaces $n$ in the above argument.
\end{proof}

\begin{remark}
{\rm
If the volume of $E$ is known, one does not need to know $N_0$
in Fact~\ref{lem:reconstruction}.
The proof of this is left to the interested reader.
\diam
}
\end{remark}

The set of Ekeland-Hofer capacities does {\it not}\, form a
generating system for symplectic capacities on $\Ell^{2n}$.
Indeed, the volume capacity $c_\vol$ is not the pointwise limit
of homogeneous monotone functions of Ekeland-Hofer capacities:

\begin{example}  \label{ex:333}
{\rm
Consider the ellipsoids $E = E(1, \dots, 1, 3^n+1)$ and $F= E(3,
\dots, 3)$ in $\Ell^{2n}$.
As is easy to see,
\begin{equation}  \label{ine:cEH}
c_k^\EH (E) \,<\, c_k^\EH (F) \qquad \text{ for all }\, k .
\end{equation}
Assume that $f_i$ is a sequence of homogeneous monotone functions
of Ekeland-Hofer capacities which converge pointwise to
$c_\vol$.
By \eqref{ine:cEH} and the monotonicity of the $f_i$ we would find
that $c_\vol (E) \le c_\vol (F)$. This is not true.
}
\end{example}

\smallskip
\begin{problem}  \label{problem:EHvolgen}
Do the Ekeland-Hofer capacities together with the volume
capacity form a generating system for symplectic capacities on $\Ell^{2n}$?
\end{problem}

If the answer to this problem is ``yes'', this is a very
difficult problem as Lemma~\ref{lem:neu} below illustrates.

\subsubsection{Ellipsoids in dimension~4}  \label{ell:dim4}

A generalized capacity on ellipsoids in dimension~4 is represented by a
function $c(a):=c\bigl(E(a,1)\bigr)$ of a single real variable $0<a\leq 1$.
This function has the following two properties.
\begin{description}
\item[(i)] The function $c(a)$ is nondecreasing.
\item[(ii)] The function $c(a)/a$ is nonincreasing.
\end{description}
The first property follows directly from the (Monotonicity) axiom. The
second property follows from (Monotonicity) and (Conformality): For
$a\leq b$, $E(b,1)\subset E \left( \frac b a a, \frac b a \right)$, hence
$c(b)\leq\frac{b}{a}c(a)$. Note that property (ii) is equivalent to the
estimate
\begin{equation}  \label{est:lip}
        \frac{c(b)-c(a)}{b-a} \,\le\, \frac{c(a)}{a}
\end{equation}
for $0<a<b$, so the function $c(a)$ is Lipschitz continuous at all
$a>0$.
We will restrict our attention to
{\em normalized}\, (generalized) capacities, so the function $c$ also satisfies
\begin{description}
\item[(iii)] $c(1)=1$.
\end{description}
An ellipsoid $E(a_1, \dots, a_n)$ embeds into $E(b_1, \dots, b_n)$ by a
{\em linear}\, symplectic embedding only if $a_i \le b_i$ for all $i$,
see \cite{HZ}. Hence for normalized capacities on the category
$\LinEll^4$ of ellipsoids with {\em linear}\, embeddings as morphisms,
properties~(i), (ii) and (iii) are the only restrictions on the function
$c(a)$.
On $\Ell^4$, nonlinear symplectic embeddings ("folding")
yield additional constraints which are still not completely
known; see~\cite{Sch} for the presently known results.

\begin{figure}[h]
 \begin{center}
  \psfrag{1}{$1$}
  \psfrag{1/2}{$\frac 12$}
  \psfrag{1/3}{$\frac 13$}
  \psfrag{1/4}{$\frac 14$}
  \psfrag{1/6}{$\frac 16$}
  \psfrag{1/8}{$\frac 18$}
  \psfrag{12}{$\frac 1 {12}$}
  \psfrag{a}{$a$}
  \psfrag{c}{$c_B(a)=a$}
  \psfrag{cv}{$c_\vol (a) = \sqrt{a}$}
  \psfrag{c2}{$\bar c_2$}
  \psfrag{l}{$l(a)$}
  \psfrag{s}{$s(a)$}
  \leavevmode\epsfbox{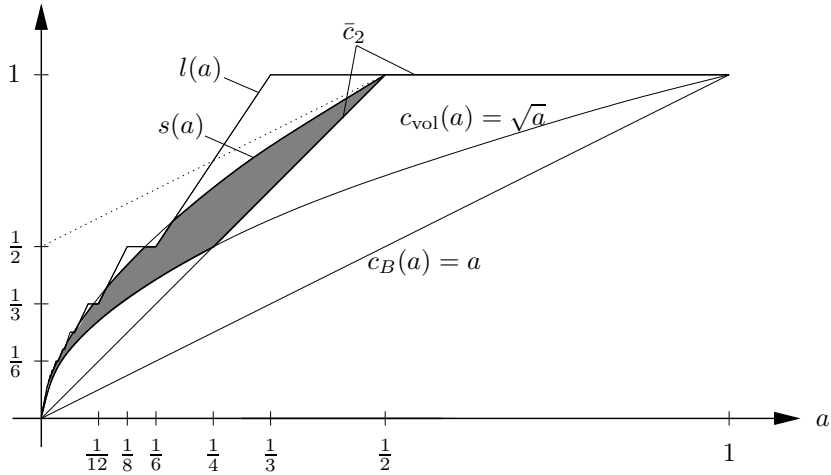}
 \end{center}
 \caption{Lower and upper bounds for $c^B(a)$.}
 \label{fi0}
\end{figure}

By Fact~\ref{fact:emb1}, the embedding capacities $c_B$ and $c^B$ are
the smallest, resp.~largest, normalized capacities on ellipsoids. By
Gromov's Nonsqueezing Theorem,
$c_B(a) = \bar c_1(a) = a$.
The function $c^B(a)$ is not completely known.
Fact~\ref{fact:emb1} applied to $\bar c_2$ yields
\[
c^B(a) =1 \,\text{ if } a \in \left[ \tfrac 12, 1 \right]
\quad \text{ and } \quad
c^B(a) \ge 2a \,\text{ if } a \in \left( 0,\tfrac 12 \right],
\]
and Fact~\ref{fact:emb1} applied to $c_\vol$ yields $c^B(a) \ge
\sqrt{a}$.
Folding constructions provide upper bounds for $c^B(a)$.
Lagrangian folding \cite{T} yields $c^B(a) \le l(a)$ where
\begin{equation*}
l(a) \,=\, \left\{
\begin{array}{ccl}
(k+1)a & \text{\rm for} & \frac{1}{k(k+1)} \le a \le
\frac{1}{(k-1)(k+1)} \\[.5em]
\frac 1 k & \text{\rm for} & \frac 1 {k(k+2)} \le a \le \frac 1{k(k+1)}
\end{array} \right.
\end{equation*}
and multiple symplectic folding \cite{Sch} yields $c^B(a) \le s(a)$
where the function $s(a)$ is as shown in Figure~\ref{fi0}.
While symplectically folding once yields $c^B(a) \le a + 1/2$ for $a \in
(0,1/2]$,
the function $s(a)$ is obtained by symplectically folding ``infinitely
many times'', and it is known that
\[
\liminf_{\eps \to 0^+} \frac{c^B \left(\frac 12 \right) -
                            c^B \left(\frac 12 -\eps \right)}{\eps}
\,\ge\, \frac 87.
\]
%

%
%
%
%
%

Let us come back to Problem~\ref{problem:EHvolgen}.
\begin{lemma}  \label{lem:neu}
If the Ekeland-Hofer capacities and the volume capacity form a
generating system for symplectic capacities on $\Ell^{2n}$,
then $c^B \left( \frac 14 \right) = \frac 12$.
\end{lemma}

We recall that $c^B \left( \frac 14 \right) = \frac 12$ means
that the ellipsoid $E(1,4)$ symplectically embeds into
$B^4(2+\eps)$ for every $\eps >0$.

{\it Proof of Lemma~\ref{lem:neu}.}
We can assume that all capacities are normalized.
By assumption, there exists a sequence $f_i$ of homogeneous and
monotone functions in the $\bar{c}_k$ and in $c_\vol$ forming
normalized capacities which pointwise converge to $c^B$.
As is easy to see,
$\bar{c}_k \left( E \left( \tfrac 14,1 \right) \right) \le
\bar{c}_k \left( B^4 \left( \tfrac 12 \right) \right)$
for all $k$,
and $c_\vol \left( E \left( \tfrac 14,1 \right) \right) =
c_\vol \left( B^4 \left( \tfrac 12 \right) \right)$.
Since the $f_i$ are monotone and converge in particular at
$E \left( \tfrac 14,1 \right)$ and $B^4 \left( \tfrac 12
\right)$ to $c^B$, we conclude that
$
c^B \left( \frac 14 \right) = c^B \left( E \left( \tfrac 14,1
\right) \right) \le c^B \left( B^4 \left( \tfrac 12 \right)
\right) = \frac 12$, which proves Lemma~\ref{lem:neu}.
\proofend

In view of Lemma~\ref{lem:neu}, the following problem is a
special case of Problem~\ref{problem:EHvolgen}.

\begin{problem}  \label{problem:12}
Is it true that $c^B \left( \frac 14 \right) = \frac 12$?
\end{problem}

The best upper bound for $c^B \left( \tfrac 14 \right)$
presently known is $s \left( \tfrac 14 \right) \approx 0.6729$.
Answering Problem~\ref{problem:12} in the affirmative means to
construct for each $\eps >0$ a symplectic embedding $E \left(
\tfrac 14,1 \right) \to B^4 \left( \frac 12 +\eps \right)$.
We do not believe that such embeddings can be constructed ``by
hand''.
A strategy for studying symplectic embeddings of \text{$4$-dimensional}
ellipsoids by algebro-geometric tools is proposed in \cite{B-ECM}.

\medskip
Our next goal is to represent the (normalized) Ekeland-Hofer
capacities as embedding capacities. First we need some preparations.

\begin{figure}[h]
 \begin{center}
  \psfrag{1}{$1$}
  \psfrag{1/2}{$\frac 12$}
  \psfrag{1/3}{$\frac 13$}
  \psfrag{2/5}{$\frac 25$}
  \psfrag{a}{$a$}
  \psfrag{co}{$c^{E(1,b)}(a)$}
  \psfrag{cu}{$c_{E(1,b)}(a)$}
  \psfrag{?}{$?$}
  \leavevmode\epsfbox{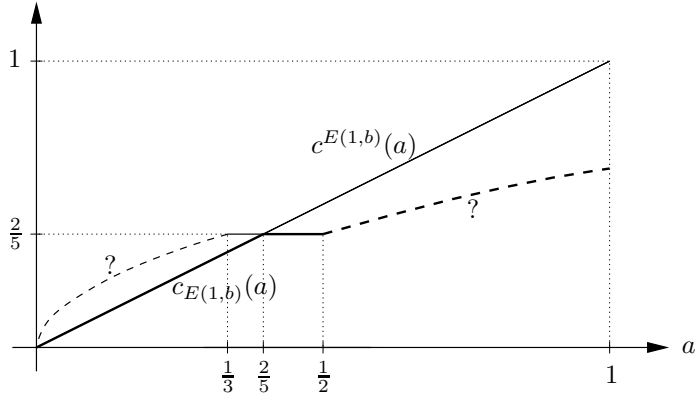}
 \end{center}
 \caption{The functions $c^{E(1,b)}(a)$ and  $c_{E(1,b)}(a)$ for $b =
 \frac 52$.}
 \label{fi2}
\end{figure}

From the above discussion of $c^B$ it is clear that capacities and
folding also yield bounds for the functions $c^{E(1,b)}$
and $c_{E(1,b)}$. We content ourselves with noting
\begin{lemma} \label{cap_ellip_4}
Let $N \in \N$ be given. Then for $N \le b \le N+1$ we have
\begin{equation}
\label{embed_to_cap}
c^{E(1,b)}(a) \,=\, \left\{
\begin{array}{cl}
\frac 1 b & \text{\rm for} \quad \frac 1 {N+1} \le a \le \frac 1 b,\\
a & \text{\rm for} \quad \frac 1 b \leq a \leq 1
\end{array} \right.
\end{equation}
and
\begin{equation}
\label{embed_from_cap}
c_{E(1,b)}(a) \,=\, \left\{
\begin{array}{cl}
a & \text{\rm for} \quad 0< a \leq \frac 1 b,\\[.5em]
\frac 1 b & \text{\rm for} \quad \frac 1 b \le a \le \frac 1 {N},
\end{array} \right.
\end{equation}
\end{lemma}
see Figure~\ref{fi2}.

\begin{remark}
Note that (\ref{embed_from_cap}) completely describes $c_{E(1,b)}$ on the
whole interval $(0,1]$ for $1 \le b \leq 2$.
\end{remark}
\begin{proof}
As both formulas are proved similarly, we only prove (\ref{embed_to_cap}).
The first Ekeland-Hofer capacity gives the lower bound
$c^{E(1,b)}(a) \ge a$ for all $a \in (0,1]$.
Note that for $a \ge \frac 1 b$ this bound is achieved by
the standard embedding, so that the second claim follows.

For $\frac 1 {N+1} \le a \le \frac 1 N$ we have $\bar c_{N+1}(E(a,1))
= 1$ and $\bar c_{N+1}(E(1,b)) = b$. Hence by Fact~\ref{fact:emb1} we
see that $c^{E(1,b)}\ge \frac 1 b$ on this interval, and this bound is
again achieved by the standard embedding. This completes
the proof of (\ref{embed_to_cap}).
\end{proof}

%
%

%
\begin{remark}
Consider the functions
\[
e^b (a) \,:=\, c^{E(1,b)} (a), \quad a \in (0,1],\, b \ge 1 .
\]
Notice that $e^1 = c^B$.
By Gromov's Nonsqueezing Theorem and monotonicity,
\[
a \,=\, c_B(a) \,=\, c^Z(a) \,\le\, e^{b}(a) \,\le\, c^B(a) ,
\quad a \in (0,1],\, b \ge 1 .
\]
Since $e^b(a) = \left( c_{E(a,1)}\bigl(E(1,b)\bigr)\right)^{-1}$ by
equation~(\ref{eq:inverse}),
we see that for each $a \in (0,1]$ the function $b \mapsto e^b(a)$ is
monotone decreasing and continuous. By~(\ref{embed_to_cap}), it
satisfies $e^b(a) = a$ for $a \ge 1/b$.
In particular, we see that the family of graphs
$\left\{ \graph \left(e^b\right) \mid 1 \le b < \infty \right\}$
fills the whole region between the graphs of $c_B$
and $c^B$, cf.~Figure~\ref{fi0}.
\diam
\end{remark}

The normalized Ekeland-Hofer capacities are represented by piecewise
linear functions $\bar c_k(a)$.
Indeed, $\bar c_1(a) =a$ for all $a \in (0,1]$, and for $k \ge 2$
the following formula follows straight from the definition
\begin{lemma}  \label{l:EH}
Setting $m:= \left[ \frac {k+1} 2 \right]$, the function $\bar c_k \colon
(0,1] \to (0,1]$ is given by
\begin{equation}
\bar c_k(a) = \left\{ \begin{array}{cl}
\frac {k+1-i} m a & \text{\rm for }\, \frac {i-1} {k+1-i} \le a
\le \frac {i}{k+1-i}  \\ [.5em]
\frac i m & \text{\rm for }\, \frac i {k+1-i} \le a \le \frac i {k-i}
\quad .
\end{array} \right.
\end{equation}
Here $i$ takes integer values between $1$ and $m$.
\end{lemma}

\begin{figure}[h]
 \begin{center}
  \psfrag{1}{$1$}
  \psfrag{1/2}{$\frac 12$}
  \psfrag{1/4}{$\frac 14$}
  \psfrag{3/4}{$\frac 34$}
  \psfrag{a}{$a$}
  \psfrag{c1}{$\bar c_1$}
  \psfrag{c2}{$\bar c_2$}
  \psfrag{c3}{$\bar c_3$}
  \psfrag{c4}{$\bar c_4$}
  \psfrag{c5}{$\bar c_5$}
  \psfrag{c6}{$\bar c_6$}
  \psfrag{cL}{$c_\infty$}
  \leavevmode\epsfbox{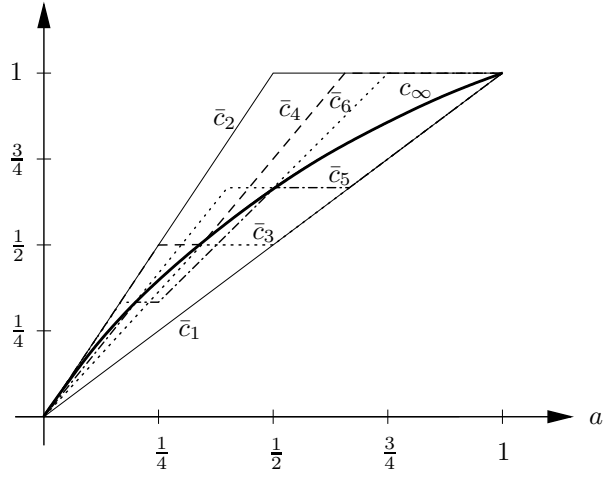}
 \end{center}
 \caption{The first six $\bar c_k$ and $c_\infty$.}
 \label{fi1}
\end{figure}

Figure~\ref{fi1} shows the first six of the
$\bar c_k$ and their limit function $c_\infty$ according to
Proposition~\ref{prop:lim_ell}.
%
%

In dimension~4, the uniform convergence $\bar c_k \to c_\infty$ is very
transparent, cf.~Figure~\ref{fi1}:
One readily checks that $\bar c_k - c_\infty \ge 0$ if $k$ is even,
in which case $\left\| \bar c_k - c_\infty \right\| = \frac{1}{k+1}$,
and that $\bar c_k - c_\infty \le 0$ if $k =2m-1$ is odd,
in which case $\left\| \bar c_k - c_\infty \right\| = \frac{m-1}{m k}$ if
$k \ge 3$.
Note that the sequences of the even
(resp.~odd) $\bar c_k$ are almost, but not quite, decreasing
(resp.~increasing). We still have
\begin{cor}  \label{l:2ml}
For all $r,s \in \N$, we have
$$
\bar c_{2rs} \le \bar c_{2r}.
$$
\end{cor}

This will be a consequence of the following characterization of
Ekeland-Hofer capacities.
\begin{lemma}  \label{l:EHext}
Fix $k\in \N$ and denote by $[a_l,b_l]$ the interval on which $\bar c_k$ has
the value $\frac l {[\frac {k+1}2]}$. Then
\begin{enumerate}[{\rm (a)}]
\item $\bar c_k \leq c$ for every capacity $c$ satisfying $\bar c_k(a_l)
\leq c(a_l)$ for all $l=1,2,\dots,[\frac {k+1}2]$.
\item $\bar c_k \geq c$ for every capacity $c$ satisfying $\bar c_k(b_l)
\geq c(b_l)$ for all $l=1,2,\dots, [\frac k 2]$ and
$$
        \lim_{a \to 0} \frac {c(a)}{a} \leq
        \frac k {\left[\frac{k+1} 2\right]}.
$$
\end{enumerate}
\end{lemma}

\begin{proof}
Formula~\eqref{est:lip} and Lemma~\ref{l:EH} show that
where a normalized Ekeland-Hofer capacity grows, it grows
with maximal slope.
In particular, going left from the left end point $a_l$ of a plateau a
normalized Ekeland-Hofer capacity drops with the fastest possible rate until
it reaches the level of the next lower plateau and then stays there, showing
the minimality.
Similarly, going right from the right end point $b_l$ of some plateau a
normalized Ekeland-Hofer capacity grows with the fastest possible rate until
it reaches the next higher level, showing the maximality.
\end{proof}

{\it Proof of Corollary~\ref{l:2ml}:}
The right end points of plateaus for $\bar c_{2r}$ are given by $b_i=\frac i
{2r-i}$. Thus we compute
$$
\bar c_{2r} \left(\frac i {2r-i} \right) = \frac i r = \frac {is}{rs} =
\bar c_{2rs} \left( \frac {is}{2rs-is} \right) = \bar c_{2rs} \left( \frac i
{2r-i} \right)
$$
and the claim follows from the characterization of $\bar c_{2r}$ by
maximality.
\proofend

Lemma~\ref{cap_ellip_4} and the piecewise linearity of the $\bar c_k$ suggest
that they may be representable as embedding capacities into a disjoint
union of finitely many ellipsoids. This is indeed the case.

\begin{prop}  \label{p:EHpresentation}
The normalized Ekeland-Hofer capacity $\bar c_k$ on $\Ell^4$ is the
capacity $c^{X_k}$ of embeddings into the disjoint union of ellipsoids
$$
X_k = Z \left( \frac m k \right) \amalg
\coprod_{j=1}^{\left[\frac k 2 \right]} E \left( \frac m {k-j}, \frac m {j}
\right),
$$
where $m= \left[ \frac {k+1} 2 \right]$.
\end{prop}
\begin{proof}
The proposition clearly holds for $k=1$.
We thus fix $k \ge 2$.
Recall from Lemma~\ref{l:EH} that $\bar{c}_k$ has $\left[\frac k 2
\right]$ plateaus, the $j^{th}$ of which has height $\frac j m$ and starts
at $a_j := \frac j {k+1-j}$ and ends at $b_j := \frac j {k-j}$.
The $j^{th}$ ellipsoid in Proposition~\ref{p:EHpresentation} is found as
follows:
In view of (\ref{embed_to_cap}) we first select an ellipsoid $E(1,b)$
so that the point $\frac 1 b$ corresponds to $b_j$.
This ellipsoid is then rescaled to achieve
the correct height $\frac j m$ of the plateau (note that by conformality,
$\alpha c^{E(\alpha, \alpha b)} = c^{E(1,b)}$ for $\alpha >0$).
We obtain the candidate ellipsoid
\[
E_j = E \left( \frac m {k-j}, \frac m j \right) .
\]
The slope of $\bar c_k$ following its $j^{th}$ plateau and the
slope of $c^{E_j}$  after its plateau both equal $\frac {k-j} m$.
The cylinder is added to achieve the correct behaviour near $a=0$.
We are thus left with showing that for each
$1 \le j \le \left[ \tfrac k 2 \right]$,
\[
\bar c_k (a) \le c^{E_j} (a) \quad \text{ for all }\, a \in (0,1] .
\]
According to Lemma~\ref{l:EHext} (a) it suffices to show that for each
$1 \le j \le \left[ \tfrac k 2 \right]$ and each
$1 \le l \le \left[ \tfrac k 2 \right]$ we have
\begin{equation}
\label{est:goal}
\bar c_k(a_l) \,=\, \frac l m \,\le\, c^{E_j}(a_l),
\end{equation}
For $l > j$, the estimate~\eqref{est:goal} follows from the fact that
$\bar c_k = c^{E_j}$ near $b_j$ and from the argument given in the proof
of Lemma~\ref{l:EHext} (a),
and for $l=j$ the estimate~\eqref{est:goal} follows from \eqref{embed_to_cap}
of Lemma~\ref{cap_ellip_4} by a direct computation.
We will deal with the other cases
$$
1 \le l< j \le \left[ \frac k 2 \right]
$$
by estimating $c^{E_j}(a_l)$ from below, using Fact~\ref{fact:emb1} with
$c=c_\vol$ and $c= \bar{c}_2$.

Fix $j$ and recall that $c_\vol(E(x,y)) = \sqrt{xy}$, so that
\begin{eqnarray*}
c^{E_j}(a_l) \geq
\frac {c_\vol(E(a_l,1))}{c_\vol \left( E \left( \frac m {k-j}, \frac m
j \right) \right)} &=&
\sqrt{\frac {l j (k-j)}{(k+1-l)m^2}}\\
&=& \frac l m \cdot \sqrt{\frac {j(k-j)}{(k+1-l)l}}
\end{eqnarray*}
gives the desired estimate \eqref{est:goal} if $j(k-j) \ge -l^2 +(k+1)l$.
Computing the roots $l_\pm$ of this quadratic inequality in $l$, we find
that this is the case if
\begin{equation*}  \label{l_minus}
l \leq l_- = \frac 1 2 \left(k+1 -  {\sqrt{1 +2k + (k-2j)^2}} \right).
\end{equation*}
%
%
%
%
Computing the normalized second Ekeland-Hofer capacity under the
assumption that $a_l \le \frac 1 2$, we find that $\bar c_2(E(a_l,1))
= 2a_l= \frac {2l}{k+1-l}$ and $\bar c_2(E_j) \leq \frac m j$, so that
\begin{eqnarray*}
c^{E_j}(a_l) \geq
\frac {\bar c_2 (E(a_l,1))}{\bar c_2 \left( E \left( \frac m {k-j},
\frac m j \right) \right)}
\geq \frac {2l}{k+1-l}\cdot \frac j m
= \frac l m \cdot \frac {2j}{k+1-l},
\end{eqnarray*}
which gives the required estimate \eqref{est:goal} if
$$
l\geq k+1-2j.
$$
Note that for $\frac 1 2 \le a_l \le 1$ we have $\bar c_2(E(a_l,1)) = 1$ and
hence
\begin{equation*}
\frac {\bar c_2(E(a_l,1))}{\bar c_2 \left( E \left( \frac m {k-j},
\frac m j \right) \right)} \ge
\frac {j} m > \frac l m
\end{equation*}
trivially, because we only consider $l < j$.

So combining the results from the two capacities, we find that the desired
estimate \eqref{est:goal} holds provided either $l \leq l_- = \frac 1 2
\left(k+1 -  {\sqrt{1 +2k + (k-2j)^2}} \right)$ or $l \geq  k+1-2j$.
As we only consider $l < j$, it suffices to verify that
$$
\min (j-1,k+1-2j) \leq \frac 1 2 \left(k+1 -  {\sqrt{1 +2k + (k-2j)^2}}
\right)
$$
for all positive integers $j$ and $k$ satisfying $1 \le j \le \left[
\frac k 2 \right]$. This indeed follows from another straightforward
computation, completing the proof of Proposition~\ref{p:EHpresentation}.
\end{proof}

Using the results above, we find a presentation of the normalized
capacity $c_\infty= \lim_{k\to \infty} \bar c_k$ on $\Ell^4$ as
embedding capacity into a countable disjoint union of
ellipsoids. Indeed, the space $X_{4r}$ appearing in the statement of
Proposition~\ref{p:EHpresentation} is obtained from $X_{2r}$ by
adding $r$ more ellipsoids. Combined with
Proposition~\ref{prop:lim_ell} this yields the presentation
$$
c_\infty \,=\, c^X \quad \text{ on }\, \Ell^4 ,
$$
where $X = \coprod_{r=1}^\infty X_{2r}$ is a disjoint union of
countably many ellipsoids.
Together with Conjecture~\ref{conj:c_L}, the following conjecture
suggests a much more efficient presentation of $c_\infty$ as an
embedding capacity.
The following result should also be proved in~\cite{CM-03b}.

\begin{conjecture}\label{conj:Lpresentation}
The restriction of the normalized Lagrangian capacity $\bar c_L$ to
$\Ell^{4}$ equals the embedding capacity $c^X$, where
$X$ is the connected subset $B(1) \cup Z(\frac 12)$ of $\R^4$.
\end{conjecture}

For the embedding capacities {\it from}\, ellipsoids, we have the
following analogue
of Proposition~\ref{p:EHpresentation}.
\begin{prop}
\label{p:EHpresentation2}
The normalized Ekeland-Hofer capacity $\bar c_k$ on $\Ell^4$ is the maximum
of finitely many capacities $c_{E_{k,j}}$ of embeddings of ellipsoids
$E_{k,j}$,
$$
\bar c_k(a) = \max \, \{ \, c_{E_{k,j}}(a)\, | \, 1 \leq j \leq m \,\}
,\quad\, a \in (0,1],
$$
where
$$
E_{k,j} = E \left( \frac m {k+1-j}, \frac m {j} \right)
$$
with $m= \left[ \frac {k+1} 2 \right]$.
\end{prop}

\begin{proof}
The ellipsoids $E_{k,j}$ are determined using \eqref{embed_from_cap} in
Lemma~\ref{cap_ellip_4}. According to Lemma~\ref{l:EHext}~(b),
this time it suffices to check that for all
$1 \leq j \le l \le \left[\frac k 2\right] $
the values of the corresponding capacities at the right end points
$b_l=\frac l {k-l}$ of plateaus of $\bar c_k$ satisfy
\begin{equation}
\label{est:goal2}
c_{E_{k,j}}(b_l) \leq \frac l m = \bar c_k(b_l).
\end{equation}
The case $l=j$ follows from \eqref{embed_from_cap} in
Lemma~\ref{cap_ellip_4} by a direct computation. For the remaining cases
$$
1 \leq j < l \leq \left[\frac k 2 \right]
$$
we use three different methods, depending on the value of $j$. If $j \leq
\frac {k-1} 3$, then Fact~\ref{fact:emb1} with $c=c_\vol$ gives
\eqref{est:goal2} by a computation similar to the one in the proof of
Proposition~\ref{p:EHpresentation}. If $j \geq \frac {k+1} 3$, then $a_j=
\frac j {k+1-j} \geq \frac 1 2$, so that \eqref{embed_from_cap} in
Lemma~\ref{cap_ellip_4} shows
that $c_{E_{k,j}}$ is constant on $[a_j,1]$, proving \eqref{est:goal2} in
this case. Finally, if $j= \frac k 3$ and $l\geq j+1$, then $\bar
c_2(E_{k,j}) = \frac {2m}{k+1-j}$ and $\bar c_2(b_l) = 1$, so that with
Fact~\ref{fact:emb1}
$$
c_{E_{k,j}}(b_l) \leq \frac {k+1-j}{2m}, 
$$
which is smaller than $\frac l m$ for the values of $j$ and $l$ we consider
here. This completes the proof of Proposition~\ref{p:EHpresentation2}.
\end{proof}

Here is the corresponding conjecture for the normalized Lagrangian
capacity.
\begin{conjecture}
The restriction of the normalized Lagrangian capacity $\bar c_L$ to
$\Ell^{2n}$ equals the embedding capacity $c_{P(1/n,\dots,1/n)}$ of
the cube of radius $1/\sqrt{n}$.
\end{conjecture}

\subsection{Polydiscs}\label{ss:poly}

\subsubsection{Arbitrary dimension}  \label{pol:arbitrary}
Again we first describe the values of the capacities in \S~\ref{s:caps}
on polydiscs.

The values of the Gromov radius $c_B$ on polydiscs are
$$
        c_B\bigl(P(a_1,\dots,a_n)\bigr)=\min\{a_1,\dots,a_n\}.
$$
As for ellipsoids, this also determines the values of $c_1^\EH$,
$c_\HZ$, $e(\cdot ,\R^{2n})$ and $c^Z$.
According to~\cite{EH-90b}, the values of Ekeland-Hofer capacities on
polydiscs are
$$
        c^\EH_k\bigl(P(a_1,\dots,a_n)\bigr) = k \pi \,\min\{a_1,\dots,a_n\}.
$$
Using Chekanov's result \cite{Che-96} that $A_\min (L) \leq e(L,\R^{2n})$
for every closed Lagrangian submanifold $L\subset \R^{2n}$, one finds the
values of the Lagrangian capacity on polydiscs to be
$$
        c_L\bigl(P(a_1,\dots,a_n)\bigr) =
        \pi \, \min\{a_1,\dots,a_n\}.
$$
Since $\vol\bigl(P(a_1,\dots,a_n)\bigr)= a_1\cdots a_n \cdot \pi^n
$ and $\vol(B^{2n})= \frac {\pi^n}{n!}$, the values of the volume
capacity on polydiscs are
$$
        c_\vol\bigl(P(a_1,\dots,a_n)\bigr)=
        \left(a_1 \cdots a_n \cdot n!\right)^{1/n}.
$$
As in the case of ellipsoids,
a (generalized) capacity $c$ on $\Pol^{2n}$ can be viewed as a function
$$
        c(a_1,\dots,a_{n-1}) := c \left( P(a_1,\dots,a_{n-1},1) \right)
$$
on the set $\{0< a_1\leq\dots\leq a_{n-1}\leq 1\}$.
Directly from the definitions and the computations above we obtain the
following
easy
analogue of Proposition~\ref{prop:lim_ell}.

\begin{prop}\label{prop:lim_pol}
As $k\to\infty$, the normalized Ekeland-Hofer capacities $\bar c_k$
converge on $\Pol^{2n}$ uniformly to the normalized Lagrangian capacity
$\bar c_L= n c_L/\pi$.
\end{prop}

Propositions \ref{prop:lim_pol} and \ref{prop:lim_ell} (together with
Conjecture~\ref{conj:c_L}) give rise to

\begin{problem}
What is the largest subcategory of $\Op^{2n}$ on which the normalized
Lagrangian capacity is the limit of the normalized Ekeland-Hofer
capacities?
\end{problem}

\subsubsection{Polydiscs in dimension~4}

Again, a normalized (generalized)
capacity on polydiscs in dimension~4 is represented by a
function $c(a):=c\bigl(P(a,1)\bigr)$ of a single real variable $0<a\le
1$, which has the properties (i), (ii), (iii).
Contrary to ellipsoids, these properties are not the only restrictions
on a normalized capacity on 4-dimensional polydiscs even if one
restricts to {\em linear}\, symplectic embeddings as morphisms. Indeed,
the linear symplectomorphism
\[
\left( z_1, z_2 \right) \,\mapsto\, \frac{1}{\sqrt{2}} \left( z_1+z_2, z_1
- z_2 \right)
\]
of $\R^4$ yields a symplectic embedding
$$
   P(a,b)\into
   P\left(\frac{a+b}{2}+\sqrt{ab},\frac{a+b}{2}+\sqrt{ab}\right)
$$
for any $a,b>0$, which implies

\begin{fact}
For any normalized capacity $c$ on $\LinPol^4$,
\[
c(a) \,\le\, \frac 12 + \frac a2 + \sqrt{a} .
\]
\end{fact}

Still, we have the following easy analogues of
Propositions~\ref{p:EHpresentation} and \ref{p:EHpresentation2}.

\begin{prop}  \label{p:EHprespol}
The normalized Ekeland-Hofer capacity $\bar c_k$ on $\Pol^4$ is the
capacity $c^{Y_k}$,
where
$$
Y_k \,=\, Z \left( \frac {[\frac {k+1} 2]} k \right)\,,
$$
as well as the capacity $c_{Y_k'}$,
where
$$
Y_k' \,=\, B \left( \frac {[\frac {k+1} 2]} k \right) .
$$
\end{prop}

\begin{cor}  \label{c:EHpres}
The identity $\bar c_k =c^{X_k}$ of Proposition~\ref{p:EHpresentation}
extends to $\Ell^4 \cup \Pol^4$.
\end{cor}

\begin{proof}
Note that $Y_k$ is the first component of the space $X_k$ of
Proposition~\ref{p:EHpresentation}.
It thus remains to show that for each of the ellipsoid components $E_j$
of $X_k$,
\begin{equation*}  \label{est:ckcj}
\bar c_k \left( P(a,1) \right)  \le\, c^{E_j} \left( P(a,1)
\right),  \quad\, a \in (0,1].
\end{equation*}
This follows at once from the observation that for each $j$ we have
$c^\EH_k \left( E_j \right) = [\frac {k+1} 2] \pi$, whereas $c^\EH_k \left(
P(a,1) \right) = ka \pi$.
\end{proof}

\begin{problem}
Does the equality $\bar c_k=c^{X_k}$ hold on a larger class of open
subsets of $\R^4$?
\end{problem}


\end{document}